\documentclass[10pt]{amsart}
\usepackage{amsmath,euscript,amssymb,amsbsy,amsfonts,latexsym,amsthm,amscd,amsxtra, xfrac, xcolor}

\usepackage[all]{xy}

\usepackage{marginnote}
\usepackage[marginparwidth=2.5cm, marginparsep=1cm]{geometry}

\newtheorem{Theorem}{Theorem}[section]
\newtheorem*{Theorem*}{Theorem}
\newtheorem{Def}[Theorem]{Definition}
\newtheorem{Prop}[Theorem]{Proposition}

\newtheorem{Lem}[Theorem]{Lemma}

\newtheorem*{MainTheorem}{Main Theorem}

\theoremstyle{definition}
\newtheorem{Rem}[Theorem]{Remark}

\newtheorem*{Pf}{Proof}
\newenvironment{Proof}{\begin{Pf} \begin{upshape}} {\end{upshape} \qed\end{Pf}}

\newcommand\beqa[1]{ \begin{eqnarray} \label{#1}}
\newcommand{\eeqa}{ \end{eqnarray} }
\newcommand{\beqano}{ \begin{eqnarray*} }
\newcommand{\eeqano}{ \end{eqnarray*} }

\newcommand{\T}{ {\mathbb T}   }
\newcommand{\N}{ {\mathbb N}   }
\newcommand{\R}{ {\mathbb R}   }
\newcommand{\Z}{ {\mathbb Z}   }
\newcommand{\C}{ {\mathbb C}   }
\renewcommand{\H}{ {\mathbb H}   }
\newcommand{\K}{ {\mathcal K}   }

\newcommand{\cg}{ {\mathfrak g}   }
\newcommand{\ccG}{ {\mathcal G}   }
\newcommand{\cG}{{\mathfrak g}}

\newcommand*{\longhookrightarrow}{\ensuremath{\lhook\joinrel\relbar\joinrel\rightarrow}}

\def\beal{\begin{aligned}}
\def\enal{\end{aligned}}
\renewcommand \a {\alpha}
\newcommand \e {\varepsilon }

\renewcommand \b  {\beta}
\newcommand \hd  {\hat{\delta}}
\newcommand \td  {\tilde{\delta}}
\newcommand \bb  {{\hat{\beta}}}

\renewcommand \d {\delta}

\newcommand \m {\mu}

\newcommand \f {\varphi}

\newcommand \g {\gamma}

\newcommand \G {\Gamma}
\newcommand \s {\sigma}

\renewcommand \l {\lambda} 

\newcommand \bL {\overline{L}}

\newcommand \bx {\bar{x}}

\newcommand \bu {\bar{u}}
\newcommand \bv {\bar{v}}
\newcommand \by {\bar{y}}
\newcommand \bH {\overline{H}}
\newcommand \bbf {\bar{f}}
\newcommand \opi {\bar{\pi}}
\newcommand \bp {\hat{p}}

\newcommand \cL {{\mathcal L}}

\newcommand \cH {{\mathcal H}}

\newcommand \cS {{\mathcal S}}

\newcommand{\calM }{{\mathfrak M}  }

\newcommand \fg {{\mathfrak g}}
\def\ie{\hbox{\it i.e.,\ }}

\newcommand \rT {{\rm T}}

\newcommand \rank {{\rm rank\ }}

\newcommand \cone {G_{\infty}}

\def\Bbb{\mathbb}

\def\R{\Bbb R}
\def\C{\Bbb C}
\def\Z{\Bbb Z}
\def\T{\Bbb T}

\def\G{\Gamma}
\def\g{\gamma}

\def\lb{\lambda}

\begin{document}

\title{On the Homogenization of the Hamilton-Jacobi Equation}
\author{Alfonso Sorrentino}
\email{sorrentino@mat.uniroma2.it}
\address{Dipartimento di Matematica, Universit\`a degli Studi di Roma ``Tor Vergata'', Rome (Italy).}
\date{\today}
\subjclass[2010]{35B27, 35B40, 35F21, 74Q10 (primary), and 20F69, 22E25, 35R03, 37J50 (secondary)}

%20F18 = Nilpotent groups
%20F65 = Geometric group theory
% 20F69 = Asymptotic properties of groups
% 22E25= Nilpotent groups
% 35 B40= Asympt behaviour of solutions of PDEs
% 35B27= Homogenization
% 74Q10 = Homogenization and oscillations in dynamical problems
% 35 F21= HJ equation
%35R01 = PDE on manifolds
% 35R03= PDE on lie groups and carnot groups
%37J50 Action min measures and orbits
\maketitle

\begin{abstract}
%Since the celebrated work by  Lions, Papanicolaou and Varadhan in 1980's, there has been a considerable attention to the  homogenization problem for  Hamilton-Jacobi equation: the interest in  this question -- besides the importance of the Hamilton-Jacobi equation in many different contexts  (classical mechanics, symplectic geometry, PDEs, etc...) -- has been in recent years boosted by the manifold connections that it shares with new prominent areas of research: Aubry-Mather theory, weak KAM theory, symplectic homogenization, just to mention a few of them.

In this article, we describe  how  the  celebrated result by Lions, Papanicolau and Varadhan on the Homogenization of Hamilton-Jacobi equation can be extended beyond the Euclidean setting. 
%We first discuss how this result and a more recent one by Contreras, Iturriaga and Siconolfi can be interpreted as a particular case of a more intrinsic and geometric approach, \ie the case of Tonelli Hamiltonians which are invariant under the action of a discrete group. In particular, we point out the leading r\^ole played by the algebraic nature of the group (more specifically, its rate of growth) in driving the homogenization process and determining the features of the limit problem. 
More specifically,  we show how to obtain a homogenization result in the case of Hamiltonians that are invariant under the action of a discrete (virtually) nilpotent group (\ie with polynomial growth), following ideas of M. Gromov \cite{Gromov} and  P. Pansu \cite{Pansu}.
\end{abstract}

%%%%%%%%%%%%%%%%%%%%%%%%%%%%%%%%%%%
\section{Introduction}

\setcounter{subsection}{1}
Since the celebrated work by  Lions, Papanicolaou and Varadhan \cite{LPV}, there has been a considerable attention to the homogenization of  Hamilton-Jacobi equation and, more generally, to homogenization problems. This increasing interest is justified by the wide spectrum of applications, particularly  to all of those models characterized by the coexistence of phenomena of different scales and diverse complexities. 
Naively speaking, one aims at describing the {\it macroscopic} behaviour and the {\it global} picture of these problems, by averaging over their {\it microscopic} oscillations and neglecting  their {\it local} features:  in a pictorial sense, the goal is to single out what remains visible to a (mathematical) observer, as her/his point of view moves farther and farther away.
One of the advantage of this large scale description is that it is expected to be  easier to study and possibly to implement (for example, numerically); at the same time, this model continues to encode much interesting information on the original (non-homogenized) problem.

\medskip
In this article, we aim to describe out  how   the  result in \cite{LPV} and the more recent one  by Contreras, Iturriaga and Siconolfi \cite{CIS} can be both interpreted as a particular case of a more intrinsic and geometric approach, which will be based on ideas of M. Gromov \cite{Gromov} and  P. Pansu \cite{Pansu}. This will allow us not only to provide a more general -- and somehow intrinsic -- homogenization result (in the same spirit as \cite{CIS}), but to shed a better light on the various features and objects that are involved in this process. We believe that this point of view has not received the right amount of attention among the PDE community.

\medskip
Since this article aims to be accessible to a wider audience, with different expertise, let us first recall  what is our main object of study (the Hamilton-Jacobi equation, subsection \ref{s1}),  what we mean by homogenization (subsection \ref{ss102}) and what our main goal is (subsection \ref{maingoal} et seq.). Expert readers can skip directly to subsection \ref{secmainres}, where the main theorem is stated.\\

\subsubsection{The Hamilton-Jacobi Equation}\label{s1}
Let us start by briefly recalling what is   the classical Hamilton-Jacobi Equation. It consists of
%\footnote{Sir William Rowan Hamilton (1805-1865) introduced this equation in his studies of geometrical optics in an arbitrary medium with varying index of refraction. Later, Carl Gustav Jacob Jacobi (1804-1851) sharpened Hamilton's formulation and clarified some mathematical issues, making significant applications to classical mechanics.}
%
%
%
a first-order nonlinear partial differential equation of the form 
$$
\partial_t u(x,t) + H(x, \partial_x u(x,t) )= 0
$$
where  $H: \R^n\times \R^n \longrightarrow \R$ is called the {\it Hamiltonian} and  $(x,t)\in \R^n \times \R$ are independent variables.%\footnote{Actually,  more general versions could be considered; for instance, one could include parameters or choose the right-hand side to be a non-vanishing function. For example: 
%$$
%\partial_t u(x, \alpha,t) + H(x, \alpha,\partial_x u(x,\alpha,t)) = K(\alpha,t),
%$$
%where $(x,\alpha,t)\in \R^n\times \R^n \times \R$  are independent variables and $K:\R^n\times\R \longrightarrow \R$ plays the role of an integrable Hamiltonian.%; this approach is particularly useful in classical mechanics, for instance to find changes of coordinates that simplify the equations.}
This equation can be easily defined on a general manifold $M$:
$$\partial_t u(x,t) + H(x, \partial_x u(x,t)) = 0
$$
where  $H: T^*M \longrightarrow \R$  is now defined on the cotangent bundle of  $M$ and  $(x,t)\in M \times \R$ continue to be independent variables.

This equation has many applications in dynamical systems and  classical mechanics; its solutions, for example, 
%provide infinite families of solutions to Hamilton's  equations of motion,
% which can be used to determine the dynamics a mechanical system or an optical system in the ray approximation. See for example \cite{Arnold, CH}.
% In some special cases -- yet very interesting --  classical solutions to the Hamilton-Jacobi equation 
are related  to the existence of Lagrangian submanifolds  that are invariant under the associated Hamiltonian flow ({\it e.g.}, KAM tori) and could be used, at least in principle, as  generators of canonical (\ie symplectic) changes of coordinates that  simplify the equations of motion and make them explicitly {\it integrable}. These special solutions -- which are known to exist only under special circumstances --  are the subject of the so-called KAM theory (named after the Kolmogorov, Arnol'd, and Moser).

Besides Hamiltonian dynamics, the Hamilton-Jacobi equation arises in many other different contexts, including: PDEs, calculus of variations, control theory, optimal mass transportation problems, conservation laws, classical limits of Schr\"odinger equation, semi-classical quantum theory, etc...  %See for example \cite{Arnold, Benton, CH, Lions, Evans}.

\subsubsection{Classical Homogenization: the Euclidean Periodic Setting}\label{ss102}
Let us start by recalling a simple model of homogenization introduced and discussed by Lions, Papanicolaou and Varadhan in \cite{LPV}. Let $H:\R^n\times\R^n \longrightarrow \R$ be a $C^2$ Hamiltonian that is $\Z^n$-periodic in the space variables $x$, strictly convex in the fiber-variables/momenta $p$ (in the $C^2$ sense, \ie $\frac{\partial^2}{\partial p^2}H(x,p)$ is positive definite everywhere) and superlinear in each fiber.\footnote{Actually, the result in \cite{LPV} is stated under less stringent hypothesis on the Hamiltonian. We restrict ourselves to  this setting, because we are interested in applying methods coming from Aubry--Mather theory and Lagrangian dynamics.} 
Observe that such a Hamiltonian can be seen as the lift of a Hamiltonian defined on $\T^n\times \R^n \simeq T^*\T^n$;  using a more modern terminology,  we would call it a {\it Tonelli Hamiltonian} (see subsection \ref{sec1.2}).

The homogenization problem is related to understanding the limit behaviour of solutions to the Hamilton-Jacobi equation, when the Hamiltonian is modified in order to have faster and faster oscillatory spatial dependence. More specifically, for small $\e>0$ one considers the rescaled equation
\begin{equation}\label{eqriscalataRn}
\left\{ \begin{array}{l}
\partial_t u^\e(x,t) + H(\frac{x}{\e},  \partial_x u^\e(x)) = 0 \qquad x\in \R^n,\; t>0\\
u^\e(x,0)= f_\e( x).\\
\end{array}\right. 
\end{equation}
It is a well-known result \cite{Ishii, CL} that if $f_\e \in C(\R^n)$ and has linear growth, then there exists a unique viscosity solution %\footnote{The concept of viscosity solution concept was introduced in the early 1980's by Pierre-Louis Lions and Michael G. Crandall as a generalization of the classical concept of solution to a partial differential equation and is based on the maximum principle for PDE's. We refer interested readers, for example, to \cite{CLvisc, Evansbook}.} 
to the above problem. In \cite{LPV} the authors studied what happens to these solutions as $\e$ goes to $0$, under the (natural) condition that in the limit $f_\e$ converge uniformly to a continuous function $\bbf$ with linear growth. In particular, they proved that the solutions 
$u^\e$ converge locally uniformly to the unique viscosity solution $\bu$ of the {\it limit} (or {\it homogenized}) {\it problem}:
\begin{equation}\label{HJRn}
\left\{ \begin{array}{l}
\partial_t \bu(x,t) + \bH(\partial_x \bu(x)) = 0 \qquad x\in \R^n,\; t>0\\
\bu(x,0)= \bbf( x),\\
\end{array}\right. 
\end{equation}
where $\bH: \R^n \longrightarrow \R$ is a convex Hamiltonian (not necessarily strictly convex) which does not depend on $x$; $\bH$  is called {\it effective} (or {\it homogenized}) Hamiltonian.
While  $\bH$ is usually highly non-differentiable, nevertheless  the solution $\bu$ of the limit problem are very easy to describe, since the characteristic curves for the limit equation  are straight lines:  for each $x\in\R^n$ and $t>0$
\begin{equation}\label{HopfRn}
\bu(x,t)=\min_{y\in\R^n} \left\{\bbf(y)+t \overline{L}\left( \frac{x-y}{t}\right)
\right\},
\end{equation}
where $\overline{L}: \R^n \longrightarrow \R$ is the Legendre-Fenchel transform of $\overline{H}$ (or {\it effective Lagrangian}),  given by
$\overline{L}(v):=\max_{p\in\R^n} \left\{ p\cdot v-\overline{H}(p) \right\}.$\\

\vspace{20 pt}

\subsection{Main goal.} \label{maingoal} The  main goal in this article is to describe how to extend this result beyond the Euclidean setting. Roughly speaking (see subsection \ref{secmainres} for a mathematically more precise description of the main results), we shall describe a homogenization  result in the case of Hamiltonians that are invariant under the action of a discrete (virtually) nilpotent group (\ie with polynomial growth), following ideas of M. Gromov \cite{Gromov} and  P. Pansu \cite{Pansu}. 
As in the above mentioned result by Lions, Papanicolaou and Varadhan  \cite{LPV}, 
we shall prove that solutions  converge, in a suitable sense, to a limit variational formula that is solution to a limit problem. \\{The term {\it Homogenization} is used by different communities with  different meanings; in this work, we refer to our result as a ``Homogenization result for the Hamilton-Jacobi equation'', in the same spirit as it is done in \cite{CIS}, \cite{Viterbo}, \cite{MVZ}, etc.}
\\

In the remaing part of this introduction, in order to present our main result, we need to describe more carefully several ingredients that are involved:
\begin{itemize}
\item the Homogenized Hamiltonian (see subsection \ref{ss103});
\item the setting to consider on a general manifold (see subsections \ref{homgenmanif} and \ref{beyondab}).
\end{itemize}

\vspace{10 pt}

\subsubsection{The Homogenized Hamiltonian}\label{ss103}
The identification of the effective Hamiltonian is definitely an important step in the homogenization procedure. In \cite{LPV}, $\bH: \R^n \longrightarrow \R$ was obtained by studying the so-called {\it cell problem} (or {\it stationary ergodic Hamilton-Jacobi equation}, see equation (\ref{cellpb})), an auxiliary equation arising in the formal expansion in $\e$ of (\ref{HJRn}). 
More specifically, for each $p\in \R^n$, there exists a unique $\l\in\R$ for which the following equation admits a viscosity (periodic) solution $v:\R^n\longrightarrow \R$:
\begin{equation}\label{cellpb}
H(x, p+\partial_x v(x)) = \l \qquad {x\in \R^n}.
\end{equation}
We denote such a value $\l$ by $\bH(p)$, thus defining $\bH:\R^n\longrightarrow \R$; see \cite[Theorem 1]{LPV}.
Observe that this equation can be thought as a nonlinear eigenvalue problem with the effective Hamiltonian $\bH(p)$ and the solution $v$ playing respectively the roles of the eigenvalue and the eigenfunction. \\

This description of the effective Hamiltonian can be generalized to the case of a Hamiltonian defined on the cotangent bundle of an arbitrary compact manifold $M$ (see for example \cite{Contreras,FM}); the previous example -- because of the periodicity assumption -- corresponds to the case
$M=\T^n$ and  $T^*M=\T^n\times \R^n$. However, a general compact manifold $M$ is not necessarily parallelizable, therefore one should pay attention to how replace the role of $p$ (\ie the argument of the effective Hamiltonian) in 
(\ref{cellpb}). 

Let us start by remarking that looking for a (smooth) periodic solution $v$ of (\ref{cellpb}) is equivalent to searching a closed $1$-form on $\T^n$ with cohomology class $p$,  whose graph is contained in an energy level of the Hamiltonian $H$ (the graph of a $1$-form corresponds in fact to a section of $T^*\T^n$). Hence, the change of $p$ can be interpreted as searching for solutions with different {cohomology classes} (recall that $H^1(\T^n;\R)\simeq \R^n$). \\

Let now $M$ be a general smooth compact manifold without boundary and let $H:T^*M \longrightarrow \R$ be a Tonelli Hamilonian, \ie 
a $C^2$ Hamiltonian, which is strictly convex and superlinear in each fiber (see subsection \ref{sec1.2} for a more precise definition). Then one can define the associated effective Hamiltonian in the following way:
\begin{itemize}
\item[-] Let $\eta$ be a closed $1$-form on $M$ and let $[\eta]\in H^1(M;\R)$ denote its cohomology class.
\item[-] There exists a unique $\lambda\in \R$  for which the following equation admits a viscosity solution $v$ (see \cite{Contreras,FM}):
\begin{equation}\label{cellpbM}
H(x, \eta(x)+d_x v(x)) = \l \qquad {x\in M}.
\end{equation}
\item[-] Clearly, this unique value $\l$ does not depend on $\eta$, but only on its cohomology class. In fact, if $\eta$ and $\eta'$ are two closed $1$-forms with $[\eta]=[\eta']$, then $\eta-\eta'$ is exact and therefore there exists 
$w: M\longrightarrow \R$ such that $\eta-\eta' = dw$, which will affect only the form of the solution and not the energy value.
\item[-] We denote this unique value of $\l$ by $\bH([\eta])$, hence defining a function $$\bH: H^1(M;\R)\longrightarrow \R$$ that will be called {\it effective Hamiltonian}.
\end{itemize}

\begin{Rem}\label{rem1.1}
It is important to recall that $H^1(M;\R)\simeq \R^{b_1(M)}$, where $b_1(M)$ denotes the first Betti number of $M$. In general, differently from what happens in the euclidean periodic case, 
there is no relation between the dimension of $M$ and $b_1(M)$, so the effective Hamiltonian is defined on a space that can have a drastically larger or smaller dimension. For example, if $M$ is a surface of genus $g\geq 0$, then 
$b_1(M)=2g$, which can be arbitrary larger than $\dim M=2$.\\
\end{Rem}

\begin{Rem}
It turns out that the effective Hamiltonian is also extremely significant from a dynamical systems point of view, particularly in the study of the associated Hamiltonian dynamics by means of variational methods (what is nowadays known as Mather and Ma\~n\'e theory), where it appears many noteworthy forms and has consequently been named in different ways by the various communities: {\it minimal average action, Mather's $\alpha$-function, Ma\~n\'e critical values}, etc... 
A brief presentation of this relation will be discussed in Appendix \ref{appMather}.\\
\end{Rem}

\subsubsection{How to Homogenize on a General Manifold}\label{homgenmanif}
We want now to address the main (and very natural) question at the core of this work:
how to generalize Lions, Papanicolaou and Varadhan's result to the case of (Tonelli) Hamiltonians defined on more general  spaces, not necessarily euclidean.\\

%We want to prove an homogenization result for Hamilton-Jacobi equation for  Tonelli Hamiltonians which are invariant under the action of a discrete group, in the spirit of \cite{LPV} and \cite{CIS}.
%In the light of the homogenization steps H1-5 singled out in subsection \ref{ss102}, 
We need to address the following questions:
\begin{itemize}
\item[Q1 -] How to rescale a general Hamilton-Jacobi equation?
\item[Q2 -] How to define the {\it homogenized} Hamiltonian $\overline{H}$? In particular, on which space should the limit problem be defined?
\item[Q3 -] How to intend (and prove) the convergence of solutions to the limit one?\\
\end{itemize}

%\vspace{10 pt}

\noindent {Q1.} Let us start by discussing the first question, which is relatively  easy to address. 
As suggested in \cite{CIS}, let us observe that in the euclidean periodic case (\ie on $\T^n$) if $u^\e$ is a solution to the rescaled Hamilton-Jacobi equation (\ref{eqriscalataRn}), then $v^{\e}(x,t):= u^\e(\e x,t)$ is a $\Z^n$-periodic solution to 
\begin{equation}\label{eqriscalataRn2}
\partial_t v^\e(x,t) + H\left(x, \frac{1}{\e} \partial_x v^\e(x,t)\right) = 0 \qquad x\in \R^n,\; t>0.
\end{equation}
This equation -- which is indeed equivalent to (\ref{eqriscalataRn}) -- can be interpreted in the following geometric way: we do not rescale the space variables, but we consider a rescaled metric on it (which affects only the momenta). Observe that  equation (\ref{eqriscalataRn2}) corresponds to Hamilton-Jacobi equation on the metric space $(\T^n, d_\e:=\e d_{\rm eucl})$, where $d_{\rm eucl}$ denotes the euclidean metric on $\T^n$.
The advantage of this formulation of the problem is that this makes sense on a general metric space, whereas rescaling the space variables is possible only on spaces equipped with dilations, or having a homogeneous structure.\\

Hence, if $X$ is a smooth connected (not necessarily compact) manifold without boundary, endowed with a complete Riemannian metric $d$
and $H: T^*X \longrightarrow \R$ is a Hamiltonian, 
then for each $\e>0$ we shall consider the following {\it rescaled Hamilton-Jacobi equation}: 
\begin{equation*}
\left\{ \begin{array}{l}
\partial_t v^\e(x,t) + H(x, \frac{1}{\e} \partial_x v^\e(x,t)) = 0 \qquad x\in X,\; t>0\\
v^\e(x,0)= f_\e( x),\\
\end{array}\right. 
\end{equation*}
where $f_\e: X \longrightarrow \R$ is a function on the rescaled metric space $X_\e:=(X, d_\e:=\e d)$ (of course, more hypothesis will be needed in order to prove a homogenization result).\\

%\vspace{10 pt}

\noindent {Q2 \& Q3}. The second and third questions are definitely less straightforward and -- as we shall see later -- they will represent the core of this article (see section \ref{secasymptgeo} for Q3 and section \ref{Hom} for Q2). \\

%\vspace{10 pt}
Actually, a more fundamental (and urgent) question should be:
\begin{itemize}
\item[Q0 -] What setting to consider?\\
\end{itemize}
%We postpone this question to a later discussion, since an answer  of this  requires a better understanding of Q1-3.\\

To the best of our knowledge, the only article in which this latter issue has been addressed before is \cite{CIS}, where the authors consider a closed manifold $M$ and discuss how to extend classical homogenization results for litfs of Tonelli Hamiltonians on the cotangent bundle of its {\it abelian cover} $\widetilde{M}$ (observe that $\R^n$ corresponds to both the universal and the abelian cover of $\T^n$). In what follows,  we shall start by recalling the reasons justifying this choice and we shall later propose a different, more general, setting that we believe better suits the nature of the problem.

\subsubsection{Homogenization on the Abelian Cover}
Let us summarize the main result  in \cite{CIS}. As we have discussed in subsections \ref{ss103}, %and \ref{ss104}, 
for Tonelli Hamiltonian $T^*M$, where $M$ is a closed manifold,  there is a natural candidate for the homogenized Hamiltonian, $\bH: H^1(M;\R) \longrightarrow \R$.
However,  there is a  crucial obstacle that one has to consider: while this limit Hamiltonian $\bH$ is defined on $H^1(M;\R)$, the rescaled solutions are functions on $M$; the problem of convergence is made even more  critical  by the fact that -- as we have pointed out in Remark \ref{rem1.1} -- these spaces have different dimensions. %(this difference might be very large!).

Inspired by a strategy that had already been exploited in the context of Mather's theory (see for example \cite{CP, Mather91}), the authors propose to solve this problem by considering the lift of the Hamiltonian to the so-called (maximal free) {\it abelian cover} of $M$, {\it i.e.}, the covering space of $M$ whose fundamental group is $\pi_1(\widetilde{M})= \ker \mathfrak{h}$, where
$\mathfrak{h}$ denotes the Hurewicz homomorphism $\pi_1(M)\longrightarrow H_1(M;\R)$, and whose group of deck transformations is isomorphic to the free part of $H_1(M;\Z)$, \ie $\Z^{b_1(M)}$ with $b_1(M):=\dim H_1(M;\R)$.

Observe that the lifted Hamiltonian $\widetilde{H}: T^* \widetilde{M} \longrightarrow \R$ satisfies a sort of  ``periodicity'' property: it is invariant under the action of the group of deck transformations. Hence, one can identify some ``privileged'' directions, which are not canceled by the rescaling process. Figuratively speaking, as one looks as the manifold from far away (\ie rescales the metric by a parameter that goes to zero), while all local properties pass out of sight, these ``homological directions'' remain clearly  identifiable. 
Then, in some sense, $\widetilde{M}$ has  a (homological) ``structure'' like $\Z^{b_1(M)}$ and consequently the $\e$-rescaled metric space
has a structure  like $\e\Z^{b_1(M)}$, which, as $\e$ goes to zero, becomes more and more similar to $\R^{b_1(M)} \simeq H_1(M;\R)$.

%This idea allowed the authors of \cite{CIS} to define a ``homological'' coordinate map $G: \widetilde{M}\longrightarrow H_1(M;\R)$ (which depends, of course, on a reference point in $\widetilde{M}$ and on the choice of a basis of $H_1(M;\R)$, like a sort of ``homological compass'') and to give a meaningful notion of convergence of $(\widetilde{M},\e \tilde{d})$ to $H_1(M;\R)$ (where $\tilde{d}$ denotes the lifted metric), as well as a notion of convergence of functions on $\widetilde{M}$ to functions defined on $H_1(M;\R)$; see \cite[Section 2]{CIS} for more details.  \\

\begin{Theorem*}[Contreras, Iturriaga and Siconolfi, \cite{CIS}]
Let  $M$ be a closed manifold, $H: T^*M \longrightarrow \R$  a Tonelli Hamiltonian
and let $\widetilde{H}:T^*\widetilde{M} \longrightarrow \R$ denote the lift of $H$ to the cotangent bundle of the maximal free abelian cover $\widetilde{M}$.
Let $f_\e : \widetilde{M} \longrightarrow \R$ and $f: H_1(M;\R)\longrightarrow  \R$ be continuous maps, with  $f$ of at most linear growth, and assume that $f_\e$ converge uniformly to $f$.
Let $v_\e: \widetilde{M}\times [0,+\infty)\longrightarrow \R $ be the viscosity solution to the problem
\begin{equation*}
\left\{ \begin{array}{l}
\partial_t v^\e(x,t) + H(x, \frac{1}{\e} \partial_x v^\e(x,t)) = 0 \qquad x\in \widetilde{M},\; t>0\\
v^\e(x,0)= f_\e( x).\\
\end{array}\right. 
\end{equation*}
Then, the family of functions $v^\e$ converges locally uniformly in $\widetilde{M}\times (0,+\infty)$ to the viscosity solution $\bv: H_1(M;\R)\times [0,+\infty) \longrightarrow \R$ of the homogenized problem:
\begin{equation*}
\left\{ \begin{array}{l}
\partial_t \bv(\bx,t) + \bH(\partial_{\bx} \bv(\bx,t)) = 0 \qquad \bx\in H_1(M;\R),\; t>0\\
\bv(\bx,0)= f(\bx),\\
\end{array}\right. 
\end{equation*}
where the effective Hamiltonian $\bH: H^1(M;\R) \longrightarrow \R$ coincides with Mather's $\alpha$-function. 

Moreover, there is a representation formula for $\bv$: 
$$
\bv(\bx,t)=\min_{\by\in H_1(M;\R)} \left\{\bbf(\by)+t \overline{L}\left( \frac{\bx-\by}{t}\right)\right\} \qquad \mbox{for}\quad \bx \in H_1(M;\R),\; t>0,
$$
where $\overline{L}: H_1(M;\R) \longrightarrow \R$ denotes  the effective Lagrangian associated to $\bH$,  \ie its Legendre-Fenchel transform (also known as Mather's $\beta$-function\footnote{See Appendix \ref{appMather}.}).\\
\end{Theorem*}

\begin{Rem}
{\it i}) The notion of convergence in the above statement must be understood in the sense introduced in \cite[Section 2]{CIS}, which is reminiscent of (pointed)  Gromov--Hausdorff convergence for (non-compact) metric spaces (see section \ref{GHconv} and \cite{BBI}). \\%Very roughly speaking, using the map $G:\widetilde{M}\longrightarrow H_1(M;\R)$, one can transpose the problem to  $H_1(M;\R)\simeq \R^{b_1(M)}$.  \\
{\it ii}) Observe that while the limit function $\bv(\cdot, t)$ is defined on $H_1(M;\R)$, the homogenized Hamiltonian $\bH$ is defined on $H^1(M;\R)$. In fact, the argument of $\bH$ is the differential of $\bu(\cdot, t)$ which is an element of the dual space $(H_1(M;\R))^* \simeq H^1(M;\R)$.\\
{\it iii}) In \cite{CIS} the authors also discuss the case in which the Hamiltonian is lifted to a non-maximal free abelian cover and obtain similar results. 
\end{Rem}

This theorem automatically raises a very natural question: {\it is there any reason why one has to consider the abelian cover of $M$ and lift the problem to this space?}
\\

First of all, it is clear that the homogenization process must take place in a non-compact manifold, otherwise it would lead to a trivial result in the limit as the rescaling parameter goes to 0.  For instance, in \cite{LPV} the authors  consider the problem not on $\T^n$, as one could think, but on $\R^n$ with $\Z^n$-periodicity conditions. 
Yet,  there are many other non-compact covers of $M$ and many possible periodicity conditions: {\it what makes the abelian cover special or preferable?}\\
%\begin{itemize}
%\item[-] 
From a technical point of view, this choice has the advantage of transforming the problem into a problem on the euclidean space $\R^{b_1(M)}\simeq H_1(M;\R)$ and makes possible to exploit the homogeneity of this space to provide a meaningful notion of convergence and extend classical results.\\ %Moreover, the ideas and techniques used in the proof, are crucially dependent on the abelianity of this group (\ie of the group of the deck transformations of the abelian cover).\\
%\item[-]
On the other hand, apparently this choice seems to be the natural one, if one wishes the homogenized Hamiltonian $\bH$ to be defined on $H^1(M;\R)$, as it is reasonable to expect in the light of what discussed in subsection \ref{ss103}. However, one could ask her/himself whether the fact that $\bH$ is defined on $H^1(M;\R)$, necessarily requires  the solution $\bv$ to  be defined on $H_1(M;\R)$.
%\end{itemize}

\subsubsection{Beyond the Abelian Case}\label{beyondab}
In order to address the above issues  and try to investigate to which extent the homogenization process can be generalized to general manifolds (hence providing a satisfactory answer to Q0), it is useful to reinterpret the results in \cite{CIS} and in \cite{LPV} in a different way. 
In particular, this different point of view presents the advantage of removing  the arbitrariness of the choice of the cover to which one lifts $H$; moreover, it is closer to the spirit of the classical result in \cite{LPV} (where no lift to a covering space is involved, since the periodicity condition  is assigned a-priori).\\

Let us start by observing that the lifted Hamiltonian $\widetilde{H}:T^*\widetilde{M}\longrightarrow \R$ is still a Tonelli Hamiltonian (it is strictly convex and superlinear in each fiber) and has the property of being invariant under the action of the group of deck transformations $\G:= {\mathcal A}{\rm ut}(\widetilde{M},M) \simeq  \Z^{b_1(M)}$ (to be more precise, it is invariant under the lifted action of $\G$ to the cotangent bundle; see subsection \ref{sec1.2}). This action enjoys many good properties (see subsection \ref{sec1.1}),  which essentially come from the fact that the quotient  $\widetilde{M}\slash\Gamma$ is a  closed manifold (in this specific case, this quotient coincides with $M$) and that the covering map is {\it regular}. \\

Hence, we believe that the most  natural setting to generalize this homogenization result is the following:
\begin{itemize}
\item Let $X$ be a smooth connected (non-compact) manifold without boundary, endowed with a complete Riemannian metric $d$; in the previous case, $X=\widetilde{M}$.
\item Let $\G$ be  a finitely generated (torsion free) group, which acts smoothly on $X$ by isometries (in other words, the metric $d$ is the lift of a metric on the quotient space $ X\slash\Gamma$); in the previous case, $\G:= {\mathcal A}{\rm ut}(\widetilde{M},M) \simeq \Z^{b_1(M)}$. 
\item Suppose that this action is free, properly discontinuous and co-compact (\ie the quotient space $ X\slash\Gamma$ is a compact manifold); in the previous case, these properties were clearly satisfied.
\item Let   $H: T^*X \longrightarrow \R$ be a Tonelli Hamiltonian. The lifted action of $\G$ to $T^*X$ is given by
\begin{eqnarray*}
\G\times T^*X  &\longrightarrow& T^*X\\
(\g, (x,p)) &\longmapsto& \g\cdot (x,p) := (\g(x),  p \circ d_{\g(x)}\g^{-1} ).
\end{eqnarray*}
We require that $H$ is invariant under this action, namely
$$
H(x,p) = H(\g(x),  p \circ d_{\g(x)}\g^{-1} ) \qquad \forall\; (x,p)\in T^*X \quad {\rm and}\quad \forall\; \g\in \G.
$$
\end{itemize}
 
\vspace{10 pt}
 
\noindent Our aim in the following sections is to provide an answer to this question: {\it Is it possible to prove in this setting a homogenization result for the associated Hamilton-Jacobi equation?} \\

%Namely: 
%\begin{itemize}
%\item Do the solutions $v^\e: X\times [0,+\infty) \longrightarrow \R$ to the rescaled Hamilton-Jacobi equation (for $\e>0$ and suitable initial data) converge to some function?
%\item  In which sense should this convergence be meant? 
%\item What is the space on which the limit function is defined? 
%\item Is the limit function a solution to a homogenized Hamilton-Jacobi equation? 
%\item What is the effective (or homogenized) Hamiltonian?\\
%\end{itemize}
%

We shall see that the answer to this question profoundly depends on the algebraic nature of the group $\Gamma$, more specifically, on its rate of growth (see subsection \ref{asymptoticcone}). 
Let us  provide first a very sketchy idea on the role of $\Gamma$ (see subsection \ref{asymptoticcone} for more details). 
When we look at a metric space $(X,d)$ from ``far away'''  (\ie we rescale the metric by a small positive parameter $\e$ and let it go to zero),  although no local property survive, yet one can still try to describe the asymptotic shape and properties of this ``limit'' space, if any limit exists (see Definition \ref{defascone}).
The key point is that these asymptotic information will be the same for  metric spaces that are ``close''  enough;  more precisely, for metric spaces that are  at finite {\it Gromov-Hausdorff distance} from each other (see section \ref{secasymptgeo}). 

In our specific case, since $\Gamma$ is acting on $X$ and the action satisfies all of the above properties, then each orbit $\Gamma\cdot x_0$ can be seen as a copy of $\Gamma$ embedded in $X$; in particular, as $\e$ goes to zero, this copy becomes denser and denser in $X$ with respect to the rescaled metric.
Hence, instead of looking at the behavior of $(X,d)$, one can consider $\G$ equipped with some suitable distance $d_\Gamma$, and study how this new metric space behaves under rescaling; this limit space -- when it exists -- is called the {\it asymptotic cone of $\Gamma$}. The existence of this limit space and its uniqueness are very subtle issues; as it turns out, a positive answer to these questions strongly depends on the algebraic nature of the group, in particular on its rate of growth, which must be at most polinomial (we shall describe these results  in subsection \ref{asymptoticcone}). %In the following we shall consider discrete free groups with {\it polinomial growth}, more specifically, torsion free {\it nilpotent groups}.

%\begin{Rem}
%The fact that the growth of $\Gamma$ must be at most polynomial should not be too suprising. For example, there is a very tight link between the rate of growth of the fundamental group of a compact manifold $M$ and the 
%rate of growth of the volume of balls on its universal cover $X$ (in this case, $\pi_1(M)$ plays the role of $\Gamma$); see for example \cite{Efremovic}  and Remark \ref{efremovic}. How the volume of balls on $X$ grows is definitely a significant piece of information in order  to study the homogenization process and to show the the convergence of the rescaled solutions.\\
%\end{Rem}
%

%In light of this, it is very  clear that the leading role in driving the homogenization process and in determining the structure of the limit space and of the homogenized equation, is played by the ``periodicity'' of the Hamiltonian (\ie the acting group $\G$), and not by the fundamental domain of this action (\ie the compact quotient $\G\backslash X$). Somehow, this also better justifies the name ``homogenization''. In particular, lifting a Tonelli Hamiltonian from a compact manifold to its abelian cover, corresponds to choosing of a specific periodic extension of the Hamiltonian, yet not unique, nor necessarily more natural than others.\\
%

%%%%%%%%%%%%%%%%%%%%%%%%%%%%%%

\subsubsection{Main Results}\label{secmainres}
Let us now describe our result on the homogenization of the Hamilton-Jacobi equation. The  setting is the one that we have described in the previous subsection and, 
as we have remarked above,  it is important to impose a condition of the growth of the discrete groups $\G$; more specifically, we shall ask that $\G$ is a discrete torsion-free {\it nilpotent} group
(see subsection \ref{secstratified} for a definition of nilpotency). 
\begin{Rem}
({\it i})  Clearly this setting embraces both the case considered in \cite{LPV} ($X=\R^n$ and $\G=\Z^n$), and the one in \cite{CIS} ($X=\widetilde{M}$ and
$\G\simeq \Z^{b_1(M)}$, where $\widetilde{M}$ is the (maximal) abelian cover of a closed manifold $M$ and $b_1(M)=\rank H_1(M;\Z)$). The case of abelian subcovers -- discussed in \cite{CIS} -- can be also treated in this way.

({\it ii}) One can extend this result to  the  case  of $\G$ being a finitely generated torsion-free group with polinomial growth (see subsection \ref{asymptoticcone} II). In fact, by a theorem of Gromov \cite{Gromov} (see also subsection \ref{asymptoticcone}), $\G$ is virtually nilpotent, \ie it contains a nilpotent torsion-free subgroup $\Gamma'$ of finite index; therefore, one can consider the action of $\Gamma'$  and apply our result to it). %(see also Remark \ref{remvirtnilp}). 
\end{Rem}

When $\G$ is nilpotent,  not only its asymptotic cone exists and is unique, but it also enjoys many interesting and useful properties.  In particular (as we shall see in section \ref{secasympnil}):
 the corresponding limit metric space is a simply connected nilpotent Lie group $G_\infty$, 
 and its Lie algebra $\cg_{\infty}$ can be equipped with a one-dimensional family of dilations, which  will be extremely useful when defining and implementing the homogenization process, somehow replacing the homogeneity  of the euclidean case. More specifically: 

 \begin{itemize}
 \item[i)] The asymptotic cone of $\Gamma$ is a simply connected nilpotent Lie group $G_{\infty}:=G_\infty(\Gamma)$, in which $\Gamma$ embeds as a co-compact lattice. Moreover, its Lie algebra $\cg_{\infty}$ is stratified (see subsection \ref{subsubseccarnotgroup}).
\item[ii)] The dimension of $G_\infty$ coincides with
$$
\sum_{k=1}^{\infty} {\rm rank}\,\left( \Gamma^{(k)}/\Gamma^{(k+1)}\right),
$$
where $\Gamma^{(k)}$ are the subgroups forming the lower central series of $\Gamma$ (see (\ref{centralseries})). Observe that if $\G$ is abelian, then this corresponds to the rank of $\G$.
\item[iii)] $G_\infty$ comes equipped with a Carnot-Carath\'eodory distance $d_\infty$ and a one-dimensional family of  dilations $\delta_t$ (dilations can be also seen as automorphisms of the algebra); see subsection \ref{subsubseccarnotgroup2}. \\
\end{itemize}

We can now state our main theorem.

\begin{MainTheorem}
Let $X$ be a smooth connected (non-compact) manifold without boundary, endowed with a complete Riemannian metric $d$. 
Let $\G$ be  a finitely generated torsion-free nilpotent group, which smoothly acts on $X$ by isometries;  suppose that this action is free, properly discontinuous and co-compact. 
% \begin{itemize}
% \item[-] There exists a simply connected nilpotent Lie group $G_{\infty}$ associated to $\Gamma$, whose Lie algebra $\cg_{\infty}$ is stratified (see subsection \ref{subsubseccarnotgroup}).
%\item[-] $G_\infty$ comes equipped with a Carnot-Carath\'eodory distance $d_\infty$ and a one-dimensional family of dilations $\delta_t$ (see subsection \ref{subsubseccarnotgroup}). 
%\item[-] The dimension of $G_\infty$ coincides with
%$$
%\sum_{k=1}^{\infty} {\rm rank}\,\left( \Gamma^{(k)}/\Gamma^{(k+1)}\right),
%$$
%where $\Gamma^{(k)}$ are the subgroups forming the lower central series of $\Gamma$ (see (\ref{centralseries})).\\
%\end{itemize}
%

Let   $H: T^*X \longrightarrow \R$ be a $\G$-{invariant  Tonelli} Hamiltonian and let  $L: T X \longrightarrow \R$ be the associated $\G$-{invariant  Tonelli} Lagrangian.

For $\e>0$, let $X_\e$ denote the rescaled metric spaces $(X,d_\e:=\e d)$ and consider the rescaled Hamilton-Jacobi equation:

\begin{equation*}
\left\{ \begin{array}{l}
\partial_t v^\e(x,t) + H(x, \frac{1}{\e} \partial_x v^\e(x,t)) = 0 \qquad x\in X_\e,\; t>0\\
v^\e(x,0)= f_\e( x),\\
\end{array}\right. \\
\end{equation*}

\noindent where $f_\e: X_\e \longrightarrow \R$ are equiLipschitz with respect to the metrics $d_\e$ and, as $\e$ goes to zero, they converge uniformly on compact sets  (in the sense of Definition \ref{defconv}) to a function $\bbf: G_\infty \longrightarrow \R$ with at most linear growth.
 Then:
 
 \begin{itemize}
\item[i)] The rescaled solutions (for $x\in X_\e$ and $T>0$)
 $$
 v^\e(x,T)=\inf  \left\{
f_\e(\g(0)) + \int_0^T L (\g(t), \e \dot{\g}(t))\,dt\; \big| \; \g\in C^1([0,T], X_\e),\; \g(T)=x \right\}\\
 $$
 converge uniformly on compact sets of $G_\infty\times (0,+\infty)$ (in the sense of Definition \ref{defconv}) to a function
  
\begin{eqnarray*}
\bar{v}:  G_\infty  \times (0,+\infty) &\longrightarrow & \R\\
(\bx,T) &\longmapsto&  \inf_{\by\in G_\infty} \left\{
\bbf(\by) +  T \overline{L}\left(\
\delta_{1/T}(\by^{-1}\bx)
\right)
\right\},
\end{eqnarray*}

\noindent where $\overline{L}: G_\infty \longrightarrow \R$ depends only on the Hamiltonian $H$ (or, equivalently, on the associated Lagrangian $L$ ). In particular, $\overline{L}$ is superlinear, \ie
$$
\forall \;A>0 \qquad \exists \;B=B(A)\geq 0\;:\qquad  
\overline{L} (\bx) \geq A\; d_{\infty}(e,\bx) - B \qquad \forall\; \bx\in G_\infty
$$
and convex, namely
$$
\overline{L}\Big(\delta_\l(\bx) \cdot \delta_{1-\l}(\by)\Big) \leq \lambda \overline{L}(\bx) + (1-\l) \overline{L}(\by)\qquad 
\forall\; \lambda\in (0,1) \quad \mbox{and}\quad \forall\; \bx, \by \in G_\infty.
$$
We shall call this function {\em Generalized Mather's $\beta$-function}.

\item[ii)] For each $\bx\in G_\infty$
$$
\overline{L}(\bx):= \inf_{\sigma \in \cH_{\bx}} \int_0^1 \b(\opi(\dot{\sigma}(s)))\,ds,
$$
where $\cH_{\bx}$ denotes the set of absolutely continuous horizontal curves $\sigma: [0,1] \longrightarrow G_\infty$ connecting $e$ to $\bx$, $\b: H_1(X\slash\Gamma; \R) \longrightarrow \R$ is Mather's $\beta$-function\footnote{See Appendix \ref{appMather}.}
associated to the Lagrangian $L$ projected on $T(X\slash\Gamma)$, and 
$\opi: \cg_\infty \longrightarrow \frac{\cg_{\infty}}{[\cg_{\infty},\cg_{\infty}]} \longhookrightarrow H_1(X\slash\Gamma; \R)$

%modulo the identification of $\frac{\cg_{\infty}}{[\cg_{\infty},\cg_{\infty}]}$ with a subspace of $H_1(\Gamma\backslash X; \R)$. 

\item[iii)]  Moreover,  $\bar{v}$ is the unique viscosity solution to the following problem:

%where $\log_{\infty}$ denotes the inverse of $\exp_{\infty}$ ($G_\infty$ is simply connected), $\opi$ is the projection $\cg_{\infty} \longrightarrow \frac{\cg_{\infty}}{[\cg_{\infty},\cg_{\infty}]}$ and $\bbb: \frac{\cg_{\infty}}{[\cg_{\infty},\cg_{\infty}]} \longrightarrow \R$ is a convex and superlinear function which depends only on the Lagrangian.
%\item  In particular,  $\bar{u}$ is the unique viscosity solution to the following problem:

%\inf_{\sigma \in \cH_{\bx}} \int_0^1 \beta(\opi(\dot{\sigma}(s)))\,ds,
%$$
%where $\cH_{x}$ denotes the set of absolutely continuous horizontal curves $\sigma: [0,1] \longrightarrow G_\infty$ connecting $e$ to $x$.

$$
\left\{ \begin{array}{ll}
\partial_t \bv (\bx,t) + \overline{H}(\nabla_\cH \bv(\bx,t)) = 0 & (\bx,t) \in G_\infty \times (0,\infty)\\
\bv(\bx,0) = \bbf(\bx) & \bx\in G_\infty,\\
\end{array}
\right.
$$

\noindent where  
$\nabla_\cH \bv(\bx,t)$ denotes the horizontal gradient (with respect to the $\bx$-component) of $\bv(\cdot, t)$ and
$\overline{H}: \left(\frac{\cg_{\infty}}{[\cg_{\infty},\cg_{\infty}]}\right)^* \longrightarrow  \R$ is the convex conjugate of $\b$ restricted to the subspace $\opi(\cg_{\infty}) \subseteq H_1(X\slash\Gamma; \R)$.\\

%\item[iii)] $\bbb$ coincides with Mather's $\beta$-function associated to the Lagrangian $L$ projected on $T(\Gamma\backslash X)$, modulo the identification of $\frac{\cg_{\infty}}{[\cg_{\infty},\cg_{\infty}]}$ with a subspace of $H_1(\Gamma\backslash X; \R)$. Similarly, $\bbb^*$ can be identified with the restriction of the corresponding Mather's $\alpha$-function (or effective Hamiltonian) to a subgroup of $H^1(\Gamma\backslash X;\R)$.\\
\end{itemize}
\end{MainTheorem}

\begin{Rem}
In the abelian case, \ie when $X$ is the maximal free abelian cover of a closed manifold $M$ and $\G$ is $\Z^{b_1(M)}$, then we recover exactly the result in \cite{CIS}. Similarly, for the abelian-subcover case.\\
\end{Rem}

%\vspace{20 pt}

\subsection{Organization of the Article}
Since this article aims to be accessible to a broad audience with different expertise, we include a presentation of some of the needed material. Namely (so to help expert readers  skip the unnecessary material): 

\smallskip
\noindent - In {\sc section} \ref{sec2} we recall some background material on group actions (subsection \ref{sec1.1}) and we introduce the notion of group-invariant Tonelli Hamiltonian and Lagrangian (subsection \ref{sec1.2}).

\smallskip
\noindent - In {\sc section} \ref{secasymptgeo} we address one of the most important issues involved in our construction and related to questions Q2\&3: the asymptotic structure and geometry of rescaled metrics spaces. After having recalled some classical material on Gromov's theory of metrics spaces (subsection \ref{GHconv}), we introduce the concept of asymptotic cone of a metric space and discuss its existence and properties in the case of a finitely generated group (subsection \ref{sec2.2}). In particular, we present Gromov and  Pansu's results on the asymptotic cone of a finitely generated nilpotent group (subsection \ref{secasympnil}). Background material on the theory of nilpotent (Lie) groups is provided.

\smallskip
\noindent - In {\sc section} \ref{Hom} we present  the homogenization procedure and prove our main result. We first start by discussing properties of  solutions to the rescaled Hamilton-Jacobi equation (subsection \ref{secrescHJe}); then, we introduce a notion of convergence for functions defined on these rescaled metric spaces, which  is  achieved by means of suitably defined rescaling maps (subsection \ref{secrescmap}). Next, a crucial step is to study the convergence of rescaled Ma\~n\'e potentials and define what we  call generalized Mather's $\beta$-function: this will play the role of the effective Lagrangian (subsection \ref{secgenbeta}). Finally, by combining all these ingredients we can prove the convergence of solutions to the rescaled problem to a solution to a well-identified limit problem, and complete the proof of the main result (subsection \ref{secproofmain}).

\smallskip
\noindent - In {\sc Appendix} \ref{appMather} we brielfy  recall some basic facts on the relation between the effective Hamiltonian and Mather-Ma\~n\'e theory. We refer to \cite{SorrentinoLectureNotes} for a more comprehensive discussion.\\

%\vspace{20 pt}

\subsection*{Acknowledgements}
I wish to express my  sincere gratitude to Gonzalo Contreras, Andrea Davini and Antonio Siconolfi for sharing with me their expertise on homogenization of Hamilton-Jacobi equation and for many helpful discussions. I would also like to thank Enrico Le Donne for his kind explanation of the results in \cite{BL}.
This work has been partially supported by the MIUR Excellence Department Project awarded to the Department of Mathematics, University of Rome Tor Vergata, CUP E83C18000100006.
%the italian national project PRIN 2012 ``{\it Critical Point Theory and
%Perturbative Methods for Nonlinear Differential Equations}'' and  by the INdAM-GNAMPA project 2013 ``{\it Tecniche simplettiche, variazionali e di viscosit\`a nell'omogeneizzazione}''.
\\

%\vspace{30 pt}

\section{Setting: Group Actions and Group-Invariant Tonelli Hamiltonians}\label{sec2}

Throughout this article  $X$ will be a smooth connected  manifold without boundary, endowed with a complete Riemannian metric. We denote by $TX$ its tangent bundle and by $T^*X$ the cotangent one. 
Moreover, we denote by $\|\cdot\|$ both the norm on $TX$ and the dual norm on $T^*X$, and by $d$ the corresponding metric on $X$. \\

We want to consider a {\it group action} on $X$ and assume that the metric is invariant under this action (in other words, the group is acting by isometries). \\ %We shall assume several conditions on this action.\\

\subsection{Group action} \label{sec1.1}

Let $(\G, \cdot)$ be a group which acts on $X$ by isometries, {\it i.e.}, there exists a group homemorphism 
\begin{eqnarray*}
\f: (\G, \cdot) &\longrightarrow& ({\rm Isom_d}(X), \circ)\\
\g &\longmapsto& \g(\cdot).
\end{eqnarray*}
Observe that this determines indeed a group action on $X$; in fact, it defines a map
\begin{eqnarray*}
\G\times X  &\longrightarrow& X\\
(\g, x) &\longmapsto& \g\cdot x := \g(x)
\end{eqnarray*}
such that
	\begin{itemize}
	\item $(\g_1\cdot \g_2)(x) = \g_1(\g_2(x))$ for all $\g_1,\g_2 \in \G$  (it follows from the fact that $\f$ is a homomorphism); 
	\item $e(x) = x$ for all $x\in X$ ($e$ denotes the identity in $\G$ and it is mapped by $\f$ into the identity in ${\rm Isom_d}(X)$, which is the identity map on $X$).\\
	\end{itemize}

In particular, the metric $d$ is $\G$-{\it invariant} (\ie invariant under this action): $d(\g(x), \g(y))=d(x,y)$ for all $x,y\in X$ and $\g\in \G$ (this is a consequence of the fact that $\G$ acts by isometries).\\

We assume that this action is:
\begin{itemize}
\item {\it Free}: if $\g(x)=x$ for some $x\in X$, then $\g=e$. %In other words, for each $x\in X$ the {\it stabilizer subgroup} ${\rm Stab}_x=\{\g\in \G:\,  \g(x)=x\}$ is trivial. In particular, this implies that $\f$ is injective.
\item {\it Properly discontinuous}: for each $x\in X$ there exists a neighborhood $U_x$ such that its $\G$-translates meet $U_x$ only for finitely many $\g\in \G$: \ie
$ \g(U_x) \cap U_x \neq \emptyset$ only for finitely many $\g\in \G$. \\
\noindent Observe that being free and properly discontinuous implies that for each $x\in X$ there exists a neighborhood of $x$,  $U_x$, such that 
$ \g(U_x) \cap U_x \neq \emptyset$  for all $\g\neq e$. In particular, the quotient $X\slash\Gamma$ is a smooth Riemannian manifold and the projection  $p: X \longrightarrow  X\slash\Gamma$ is a covering map.
\item {\it Co-compact}: the quotient space $ X\slash\Gamma$ is compact (with respect to the quotient topology).\\
\end{itemize}

\subsubsection{Examples.}\label{subsecexamples} Consider a compact manifold $M$ and let $p: \widehat{M} \longrightarrow M$ be a  {\it regular} covering map  with $\widehat{M}$ being a manifold; we recall that a covering map $p: \widehat{M} \longrightarrow M$ is said to be  {\it regular}
(also called {\it normal} or {\it Galois}), if $p_*\pi_1(\widehat{M}, \hat{x}_0)$ is a normal subgroup of $\pi_1(M, x_0)$. %In particular, this implies that the group of deck transformations (or automorphisms) of this cover $p$ form a group with respect to the composition: $({\mathcal A}{\rm ut}(p), \circ)$.
This implies that  the action of the group of deck transformations (or automorphisms) ${\mathcal A}{\rm ut}(p)$  is free and  transitive on all fibers. In particular, ${\mathcal A}{\rm ut}(p)$ is isomorphic to a subgroup of $\pi_1(M)$; more precisely, 
$$
{\mathcal A}{\rm ut}(p) \simeq \frac{\pi_1(M, x_0)}{p_*(\pi_1(\widehat{M},\hat{x}_0))}.
$$
\begin{itemize}
\item {\it Universal cover:} We take $X=\widetilde{M}$, the universal cover of $M$, and $\Gamma={\mathcal A}{\rm ut}(p) \simeq \pi_1(M)$.
\item {\it Abelian cover:} We take $X=\widehat{M}$, the abelian cover of $M$, {\it i.e.}, the covering space of $M$ whose fundamental group is $\pi_1(\widehat{M})= \ker \mathfrak{h}$, where
$\mathfrak{h}$ denotes the Hurewicz homomorphism $\pi_1(M)\longrightarrow H_1(M;\Z)$. In this case the group of deck transformations is given by ${\mathcal A}{\rm ut} = {\rm Im}\, \mathfrak{h}\simeq [\pi_1(M),\pi_1(M)]$, that is the commutator subgroup; in particular $\Gamma$ is an abelian group: 
$$
\Gamma \simeq \frac{\pi_1(M)}{[\pi_1(M),\pi_1(M)]} \simeq H_1(M;\Z).\\
$$  

\begin{Rem}
Hereafter we shall refer to the {\it maximal free abelian cover} as the covering space $X=\widehat{M}$ whose fundamental group is $\pi_1(\widehat{M})= \ker \mathfrak{h}$ and whose group of Deck transformation is isomorphic to the free part of 
$H_1(M;\Z)$, \ie it is isomorphic to $\Z^{b_1(M)}$ where $b_1(M)$ denotes the first Betti number of $M$.
\end{Rem}
\end{itemize}
\vspace{10 pt}

\subsubsection{Metrics on $\Gamma$}\label{metrics} One can define several different metrics on $\Gamma$:
\begin{itemize}
\item[{ I.}] 
With the above conditions on the action, it follows that every $\G$-invariant metric $d$ on $X$ and every point $x\in X$ determine a left-invariant metric on $\G$, called an {\it orbit metric}:
$$
d_{\G,x}(\g_1,\g_2) = d(\g_1(x), \g_2(x)).
$$
In other words, one identifies the group $\G$ with the orbit $\G\cdot x$: the metric $d_{\G,x}$ is nothing else than the metric $d$ restricted to $\G\cdot x$. \\
Note that if $X$ is a length space with the above $\G$-action and we consider the projection map $p: X \longrightarrow X\slash\Gamma$, then the group $\G$ with the metric $d_{\G,x}$ is isometric to $p^{-1}(y)$ with the metric induced from $X$. In particular, $p^{-1}(y)$ is a separated net\footnote{Recall that a set $S \subset X$ is called a {\it net} in $X$ if the Hausdorff distance between $X$ and $S$ is finite (see Definition \ref{Hdistance}).
$S$ is {\it separated net} if it is a net and  there exists $\e>0$ such that $d(x_1,x_2)\geq \e$ for all $x_1,x_2\in S$; in particular, $S$ is also said $\e$-{\it separated}.
}
in $X$ (where $y=p(x)$). \\ 

\item[{ II.}] If $\Gamma$ is finitely generated, one can introduce the notion of {\it word metric}. Let 
$S=\{s_1, \ldots, s_k\}$ be  a symmetric generating set (symmetric means that if $s \in S$ then also its inverse $s^{-1}$ belongs to $S$).  For each $\g\in \G$ we define the algebraic norm
$$
\|\g\|_S := \min \{n\in \N:\; \gamma \in S^n\},
$$
that is  the smallest  $m\in \N$ such that  $\g= s_{i_1}\cdot \ldots \cdot s_{i_m}$, with $s_{i_j}$ elements of $S$.
The {\it word distance} between $\g_1$ and $\g_2$ is given by $\rho_S(\g_1,\g_2):= \|\g_1^{-1}\gamma_2\|_S$. Clearly it is a metric and it is also left-invariant.\\

\end{itemize}

It is possible to relate these different kinds of metrics  (see \cite[Theorem 8.3.19]{BBI}). 

\begin{Prop}\label{prop1}
Let $\G$ be a finitely generated group and $d_{\G,x}$ be an orbit metric of a free, co-compact action of $\G$ by isometries on a length space $X$. Then, $d_{\G,x}$ is bi-Lipschitz equivalent to a word metric. Since all word metrics on $\G$ are bi-Lipschitz equivalent to one another, then all such orbit metrics and word metrics on $\G$ are bi-Lipschitz equivalent to one another.\\
\end{Prop}

\subsection{Group-Invariant Tonelli Hamiltonians and Lagrangians} \label{sec1.2}

First of all, observe that the above action of $\G$ on $X$ by isometries naturally extends to an action on $T^*X$  and $TX$, by means of the differentials of these maps:

\begin{eqnarray*}
\G\times T^*X  &\longrightarrow& T^*X\\
(\g, (x,p)) &\longmapsto& \g\cdot (x,p) := (\g(x),  p \circ d_{\g(x)}\g^{-1} )
\end{eqnarray*}

and

\begin{eqnarray*}
\G\times T X  &\longrightarrow& T X\\
(\g, (x,v)) &\longmapsto& \g\cdot (x,v) := (\g(x),  d_{x}\g(v) ).\\
\end{eqnarray*}

%\vspace{20 pt}

Now,  let $H: T^*X \longrightarrow \R$ be a $\G$-{\it invariant  Tonelli} Hamiltonian. More specifically:\\

\begin{itemize}
\item[(i)] $H\in C^2(T^*X)$;
\item[(ii)] $H$ is $C^2$-strictly convex in each fibre, {\it i.e.}, $\frac{\partial^2 H}{\partial p^2}(x,p)$ is strictly positive definite  for every $(x,p)\in T^*X$;
\item[(iii)] $H$ is uniformly superlinear in each fibre:  for every $A\geq 0$, there exists $B(A)\in \R$ such that
$$
H(x,p) \geq A\|p\|_x - B(A) \qquad \forall \; (x,p)\in T^*X;
$$
\item[(iv)] $H$ is $\G$-invariant, {\it i.e.}, $H$ is invariant under the action of $\G$ on $T^*X$:
$$
H(\g(x), p \circ d_{\g(x)}\g^{-1}) = H(x,p) \qquad \forall\, (x,p) \in T^*X \; {\rm and}\; \forall\, \g\in \G.
$$
%(or equivalently
%$H(\g(x), p) = H(x,p \circ d_x\g)$ for every  $(x,p) \in T^*X$ and $ \g\in \G$).\\
\end{itemize}

\noindent In other word, $H$ is the lift of a Tonelli Hamiltonian on $T^*( X\slash\Gamma)$. 
Since $X\slash\Gamma$ is  a compact manifold, then it follows from the above conditions that the Hamiltonian flow of $H$ is complete.\\

Let us now consider its associated Lagrangian $L: TX \longrightarrow \R$, defined as:
$$
L(x,v) = \sup_{p\in T^*_xX} \left( \langle p, v \rangle - H(x,p) \right).
$$
One can easily check that $L$ is still of Tonelli type ({\it i.e.}, it satisfies (i)--(iii)). 
Moreover, it follows from  (iv) that $L$ is also $\G$-invariant . In fact, for all $(x,v) \in TX$ and $\g\in \G$ we have: 
\begin{eqnarray*}
L(\g(x), d_x\g(v)) &=&
\sup_{p\in T^*_{\g(x)}X} \left( \langle p, d_x\g(v) \rangle - H(\g(x),p) \right) 
\;=\;
\sup_{p\in T^*_{\g(x)}X} \left( \langle p \circ d_x\g, v \rangle - H(x,p \circ d_x\g) \right) \\
&=&
\sup_{\tilde{p}\in T^*_{x}X} \left( \langle \tilde{p} , v \rangle - H(x,\tilde{p}) \right) 
= L(x,v).\\
\end{eqnarray*}

%\vspace{20 pt}

%%%%%%%%%%%%%%%%%%%%%%%%%%%%%%%%%%%%%%%%%%%%
%%%%%%%%%%%%%%%%%%%%%%%%%%%%%%%%%%%%%%%%%%%%

\section{Asymptotic Geometry of the Rescaled Spaces}\label{secasymptgeo}

In this section we would like to provide a concise presentation of some classical material related to our discussion of  questions Q2\&3.  We refer interested readers to  \cite{BBI}   for a more comprehensive discussion of these (and many other related) topics.\\

\subsection{Gromov--Hausdorff Convergence of Metric Spaces} \label{GHconv}
Let us start by recalling the definition of {\it distance} between metric spaces and use it to introduce a meaningful notion of convergence.\\

\begin{Def}[{\bf Hausdorff distance}] \label{Hdistance} Let  $(X,d)$  be a metric space and let $A$ and $B$ be two subsets of $X$. The {\it Hausdorff distance} between them, denoted $d_{H}(A,B)$, is defined by 
$$
d_H(A,B)= \inf \left\{ r>0: \; A\subset U_r(B)\; {\it and}\; B\subset U_r(A) \right\},
$$
where $U_r(A)$ denotes the $r$-neighborhood of $A$, {\it i.e.}, the set of points $x$ such that $d(x,A):=\inf\{d(x,y):\; y\in A\} < r$ {\rm (}similarly for $U_r(B)${\rm )}.\\
\end{Def}

\begin{Def}[{\bf Gromov--Hausdorff distance}] Let $X_1=(X_1,d_1)$ and $X_2=(X_2,d_2)$ be two metric spaces. The {\it Gromov--Hausdorff distance} between them, denoted $d_{GH}(X_1,X_2)$, is defined by the following relation. For any $r>0$, $d_{GH}(X_1,X_2)<r$ if and only if there exist a metric space $Z=(Z,d_Z)$ and subspaces $X_1'$ and $X_2'$ of it which are isometric, respectively, to $X_1$ and $X_2$ and  such that $d_H(X_1',X_2')<r$. In other words, $d_{GH}(X_1,X_2)$ is the infimum of positive $r$ for which the above $Z$, $X_1'$ and $X_2'$ exist.\\
\end{Def}

\begin{Rem} \label{remark1} Some properties of the Gromov--Hausdorff distance:\\

\noindent(i) $d_{GH}$ satisfies the triangle inequality: 
$$
d_{GH} (X_1,X_2) \leq d_{GH} (X_1,X_3) + d_{GH} (X_3,X_2)
$$
for any metric spaces $X_1$, $X_2$ and $X_3$.\\

\noindent (ii) It is easy to check, using the definition, that the Gromov--Hausdorff distance between isometric spaces is zero. In particular, in the case of compact metric spaces the converse is true: if $X_1$ and $X_2$ are compact metric spaces such that $d_{GH}(X_1, X_2)=0$, then they are isometric. Hence, $d_{GH}$ defines a finite metric on the space of isometry classes of compact metric spaces (see \cite[Theorem 7.3.30]{BBI}). \\
One can actually prove something more: if $X_1$ is a compact metric space and $X_2$ is a complete metric space such that $d_{GH}(X_1,X_2)=0$, then $X_1$ and $X_2$ are isometric.\\

\noindent (iii) Let us first recall the definition of $\e$-isometry. Let $(X_1,d_1)$ and $(X_2,d_2)$ be two metric spaces and let $\e>0$; $f: X_1 \longrightarrow X_2$ is an $\e$-isometry if:
\begin{itemize}
\item[1.] ${\rm dis}(f):= \sup_{x,x' \in X_1} \left| d_2(f(x),f(x')) - d_1(x,x') \right| \leq \e$   (${\rm dis}(f)$ is called {\it distortion of $f$});
\item[2.] $f(X_1)$ is a $\e$-net in $X_2$, {\it i.e.}, for every $x_2\in X_2$ we have $d_2(x_2, f(X_1)) \leq \e$.\\
\end{itemize}

\noindent We can now state a sort of  generalisation of property (iii) above (see \cite[Corollary 7.3.27]{BBI}). Let $(X_1,d_1)$ and $(X_2,d_2)$ be two metric spaces. The following result is true: 
\begin{itemize}
\item if $d_{GH}(X_1,X_2) < \e$, then there exists a $2\e$-isometry from $X_1$ to $X_2$;
\item if there exists an $\e$-isometry from $X_1$ to $X_2$, then $d_{GH}(X_1,X_2)<2\e$.\\
\end{itemize}

\noindent(iv) It follows from the definition, that if $X_2$ is an $\e$-net in a metric space $X_1=(X_1,d_1)$, then $d_{GH}(X_1,X_2)\leq \e$, where the metric on $X_2$ is the metric induced from $d_1$. In fact, it is sufficent to take $Z=X_1$, $X'_1= X_1$ and $X'_2=X_2$.
In particular, let $X$ be a length space and $\G$ be a group acting on it (we assume, as usual, the action to be free, properly discontinuous and co-compact). Then, for any $y\in  X\slash\Gamma$, $p^{-1}(y)$ is a separated net in $X$ and the group $\G$ with the orbit metric $d_\G$ is isometric to $p^{-1}(y)$ with the metric induced from $X$. In particular,  the Gromov--Hausdorff distance between $X$ and the group is  finite.\\ 
\end{Rem}

%\vspace{10 pt}

We would like to introduce a notion of {\it convergence of metric spaces}. One could easily consider the notion given by $d_{GH}$. While this works well for compact metric spaces, for non-compact  ones a slighlty more general notion is needed\footnote{The situation is similar to  what happens with the uniform convergence of functions on a fixed domain. If the domain is compact, then uniform convergence is a widely used notion; however, it becomes very restrictive once non-compact domains come into questions. For, one introduces the notion of  uniform convergence on compact sets: a sequence of functions converges if it converges uniformly on every compact subset of the domain.}.
Roughly speaking,  a sequence $\{(X_n,d_n)\}$ of metric spaces converges to a space $(X,d)$ if for every $r>0$ the balls of radius $r$ in $X_n$ centered at some fixed points converge (as compact metric spaces) to a ball of radius $r$ in $X$.\\

Let us state this convergence more precisely. First of all, recall that a {\it pointed metric space} is a triple $(X,d,x)$, where $(X,d)$ is a metric space and $x$ a point in $X$.\\

\begin{Def}[{\bf Pointed Gromov-Hausdorff convergence}]
A sequence of pointed metric spaces $\{(X_n,d_n,x_n)\}$ converges in the Gromov-Hausdorff sense to a pointed metric space $(X, x)$ --  which we shall denote
$
(X_n,d_n,x_n) \stackrel{{\rm GH}}{\longrightarrow} (X,d,x)
$
--  if the following holds. For every $r>0$ and $\e>0$ there exists 
$n_0\in\N$ such that for every $n>n_0$ there is a (not necessarily continuous) map $f_n: B_r(x_n) \longrightarrow X$ such that the following hold:
\begin{itemize}
\item[(1)] $f_n(x_n)=x$;
\item[(2)] ${\rm dis}(f_n)<\e$  {\rm (recall the definition of {\em distortion} in Remark \ref{remark1} {\rm (iii)});}
\item[(3)] the $\e$-neighborhood of the set $f_n(B_r(x_n))$ contains the ball $B_{r-\e}( x).$\\
\end{itemize}
\end{Def}

\begin{Rem}
\noindent(i) For compact metric spaces this convergence is equivalent to the ordinary Gromov--Hausdorff convergence.\\

\noindent(ii) Requirements (1) and (2)  in the above definition imply that the image $f_n(B_r(x_n))$ is contained in the ball of radius $r+\e$ centered at $x$. In particular, this and requirement (3)  imply (see \cite[Corollary 7.3.28]{BBI}) that the ball $B_r(x_n)$ in $X_n$ lies within the Gromov--Hausdorff distance of order $\e$ from a subset of $X$ between the balls of radii $r-\e$ and $r+\e$ centered at $x$ (here ``between'' means that the sets contains one ball and is contained in the other).\\

\noindent(iii) Obviously, if a sequence of pointed metric spaces converges to a pointed metric space $(X,d,x)$, then it also converges to its completion. Hence, we shall always consider complete metric spaces as Gromov--Hausdorff limits. Then, a Gromov--Hausdorff limit of pointed spaces is essentially unique (see \cite[Theorem 8.1.7]{BBI}):  let $(X,d,x)$ and $(X',d',x')$  be two (complete) Gromov--Hausdorff 
limits of a sequence $\{(X_n,d_n,x_n)\}$, and assume that $X$ is boundedly compact ({\it i.e.}, all closed and bounded sets are compact).
Then, there exists an isometry $f:X \longrightarrow X'$ with $f( x)=x'.$\\

\noindent(iv) Let $(X_n,d_n,x_n)\stackrel{{\rm GH}}{\longrightarrow} (X,d,x)$, where $X_n$ are length spaces  ({\it i.e.}, the metric is obtained from a {\it length structure}) and $X$ is complete, then $X$ is also a length space.\\

\noindent(v) If $X$ is a length space, property (ii)  can be made more precise.  In fact, one can show that for every $r>0$ the $r$-balls in $X_n$ centered at $p_n$ converge (with respect to the Gromov--Hausdorff distance) to the $r$-ball in $X$ centered at $x$.\\
\end{Rem}

%%%%%%%%%%%%%%%%%%

\subsection{Large Scale Geometry  and Asymptotic Cones}\label{sec2.2}

In this subsection we want to study the large scale geometry of a metric space. Roughly speaking, we shall look at a metric space from ``far away'' and try to describe its asymptotic  shape and properties. As we shall see, no local properties will survive, while  asymptotic properties will be the same for spaces at finite Gromov--Hausdorff distance from each other.\\

Let us recall that for a metric space $X=(X,d)$ and $\e>0$ one can consider a {\it rescaled} metric space $X_\e=(X, \e d)$, {\it i.e.}, the same set of points equipped with a rescaled metric. Similarly, for a pointed metric space.
In particular, a pointed metric space $(X,d,x)$ is called a {\it cone} if it is invariant under rescaling, {\it i.e.}, for every $\e>0$ we have that $(X, \e d, x)$ is isometric to $(X, d, x)$ as a pointed space.
\\

\begin{Def}[{\bf Asymptotic Cone}] \label{defascone} Let $(X,d)$ be a (compactly bounded) metric space and $x\in X$. A Gromov--Hausdorff limit of pointed spaces $(X, \e d, x)$ as $\e\rightarrow 0$, if any exists, is called a {\it Gromov--Hausdorff asymptotic cone of} $X$ at infinity.\\
\end{Def}

\begin{Rem}
\noindent (i) The asymptotic cone is a cone and it does not depend on the choice of the reference point $x$ (see \cite[Proposition 8.2.8]{BBI}).\\

\noindent (ii) Let $(X_1,d_1)$ and $(X_2,d_2)$ be two metric spaces such that $d_{GH}(X_1,X_2)<\infty$. Then, if $X_1$ has an asymptotic cone, then $Y$ has one too, and the two cones are isometric. In particular, if a metric space $X_1$ lies within finite Gromov--Hausdorff distance from some cone $Y$, then $Y$ is an asymptotic cone of $X$.\\

\noindent (iii) Not all spaces have an asymptotic cone. For example, the hyperbolic plane $\H^2$ with Poincar\'e's metric has no asymptotic cone; the reason, in plain words,  is that its metric balls grow too fast when the radius goes to infinity (see \cite[Exercise 8.2.13]{BBI}).\\
\end{Rem}

We have already pointed out that if  $X$ is a length space and $\G$ is a group acting on it (we assume, as usual, the action to be free, properly discontinuous and co-compact), then the Gromov--Hausdorff distance between $X$ and the group (with any orbit or word metric) is finite. Hence, for large scale considerations one can replace $X$ by $\G$. \\
This can be stated more precisely.
First, let us recall that two metric spaces  $(X_1,d_1)$ and $(X_2,d_2)$  are {\it quasi-isometric} if there exists a map $f:X_1 \longrightarrow X_2$ such that
$$
\lambda^{-1} d_1(x,y) - C \leq  d_2(f(x_1), f(x_1')) \leq \lambda\, d_1(x,y) + C
$$
for some  constants $C\geq 0$ and $\lambda\geq 1$ and for every $x_1,x_1' \in X_1$.\\
The above considerations  imply the following (see also \cite[Corollary 8.3.20]{BBI} and Proposition \ref{prop1}).\\

\begin{Prop}\label{prop2}
All length spaces $X$ admitting a free, properly discontinuous, co-compact action of a given group $\G$ are quasi isometric to one another, and are quasi-isometric to the group $\G$ equipped with any (word or orbit) metric.\\
\end{Prop} 

\noindent In particular, being quasi-isometric implies that the spaces have finite Gromov--Hausdorff distance; hence, using what remarked above,  the asymptotic cone of $X$ is isometric to the asymptotic cone of the group $\G$. In the following subsection we shall discuss the asymptotic cone of a group.\\

\subsubsection{Asymptotic Cone of a Finitely Generated Group.}\label{asymptoticcone}
%This notion has been introduced by Gromov in \cite{Gromov} to prove his famous theorem on groups with polynomial growth. 
Let $\G$ be a finitely generated group equipped with a (word or orbit) metric $d$ (see subsection \ref{metrics}). Consider the sequence of metric spaces $\G_\e=(\G, d_\e)$ where $d_\e(\g_1,\g_2) = \e\, d(\g_1,\g_2)$. 
An {\it asymptotic cone of $\G$}, which we shall denote  $\cone(\G)$ (or simply $G_{\infty}$, if there is no ambiguity on $\Gamma$), corresponds to a Gromov-Hausdorff limit of the metric spaces $\{G_\e\}$, if any limit exists. Because of Proposition \ref{prop2}, the choice of the metric does not play any significant  role in its definition (up to isometry).\\

\noindent This raises the following {questions:} {when does an asymptotic cone exist? If it exists, is it unique?\\}

\noindent These questions have different answers according to the algebraic nature of the group.\\

\noindent {\bf I.} 
Let us start with the easier case of {\bf finitely generated abelian groups}. If $\G$ is a finitely generated abelian group, then $\G$ can be decomposed into the direct sum $\Z^k\oplus \G_0$  where $k\geq 0$ is the {\it rank} of $\G$ and $\G_0$ is some finite group (consisting of elements of finite order), called the {\it torsion subgroup}. Since  $\G$ is at finite Hausdorff distance from its $\Z^k$ component, then the asymptotic cone of $\G$ is the same as the asymptotic cone of $\Z^k$.  Let $V$ be the ambient vector space, obtained from $\G$ by tensor multiplication, $G_{\infty}=\G\otimes \R$ 
(in other words, one can assume without any loss of generality that $\G\simeq \Z^k$ and $G_{\infty}\simeq \R^k$). Then: %as the integer points of $\R^k$, one can easily deduce that:
\begin{itemize}
\item[-] for any $\g\in \G$, the limit $\lim_{n\rightarrow +\infty} \frac{d(0, n\g)}{n}$ exists (see \cite[Proposition 8.5.1]{BBI});
\item[-] there exists a unique norm $\|\cdot \|_\infty$ on $G_{\infty}$, such that $\|\g\|_\infty=\lim_{n\rightarrow +\infty} \frac{d(0, n\g)}{n}$ for all $\g\in \G$ (the so-called {\it stable norm } of $\G$); 
\item[-] the asymptotic cone of $(\G,d)$ is isometric to $(G_{\infty},d_{_\infty})$, where $d_{_\infty}$ is the distance induced by ${\|\cdot\|_\infty}$.\\
\end{itemize}

\noindent{\bf II.}
 More generally, let us see what happens in the case of {\bf finitely generated groups with polinomial growth.} Let $\G$ be a finitely generated group with a word metric $d$ and let us denote by $e$ its identity element. We  say that $\G$ has {\it polinomial growth} if there exist $C>0$ and $K>0$ such that
for each $r>0$ we have
$$
\sharp \{\g\in \G: \; d(e,\g)\leq r \}  \leq C r^K.
$$
It is possible to check that this notion is well-defined and that the definition does not depend on the choice of the word metric \cite{Gromov}.\\

\noindent {\bf Examples:}\\ 
i) Finitely generated abelian groups have clearly polynomial growth. It is an easy exercise to check that the growth rate equals the rank of the group).\\
ii) Finite extensions of groups with polynomial growth, have polynomial growth.\\
iii) A finitely generated group with a nilpotent subgroup of finite index, has polynomial growth (this was proved by  Wolf in \cite{Wolf}; see also Bass' result in \cite{Bass}). The growth rate corresponds to what is called the {\it homogeneous dimension of the group} (see (\ref{homdim})). \\

\noindent Gromov in \cite{Gromov} proved the following converse result: \\

\noindent{\bf Theorem} (Gromov \cite{Gromov}).
{\it If a finitely generated group $\G$  has polynomial growth, then it contains a nilpotent subgroup of finite index}.\\

In the following we shall refer to these groups as {\it virtually nilpotent}, \ie they are finitely generated and  contain a nilpotent subgroup of finite index (or equivalently, they have polinomial growth).
In this case,  Gromov \cite{Gromov} and Pansu \cite{Pansu}  provided a very precise description of the asymptotic cone and of its properties; we shall discuss it in more details in the next subsection.
\\

\begin{Rem}\label{efremovic}
Recall that the fundamental group of a compact manifold $X$ has polynomial growth if and only if  the universal covering $Y$ of $X$ has {\it polynomial growth}\footnote{
Consider a Riemannian manifold $Y$ and for $y\in Y$ denote by ${\rm Vol}_y( r)$ the volume of the ball of radius $r$ around $y$. The {\it growth } of $Y$ is defined as the asymptotic behaviour of ${\rm Vol}_y( r)$ as $r$ goes to infinity. Efremovic \cite{Efremovic} pointed out that the growth of a manifold $Y$, which covers a compact manifold $X$, only depends on the fundamental groups $\pi_1(X)$, $\pi_1(Y)$ and the inclusion $\pi(Y) \subset \pi_1(X)$.} \cite[Section 2]{Gromov}. Observe that most (complete, non-compact) manifolds have exponential growth, but there are some interesting examples of manifolds of polynomial growth:
\begin{itemize}
\item[a)] complete manifolds of non-negative Ricci curvature; 
\item[b)] real algebraic submanifolds in $\R^n$;
\item[c)] nilpotent Lie groups with left invariant metrics; 
\item[d)] leaves of Anosov foliations.\\
\end{itemize}
%\end{Rem}
%
%\begin{Rem}\label{amenability}
%ii) Observe that if $\G$ has polinomial growth, then it is {\it amenable}. In fact, non-amenable groups have exponential growth (see \cite[Lemma A.12]{Kleiner}). Let us recall that a group $\G$ is {\it amenable} if  (there are many equivalent definitions) there exists a right invariant mean on $\ell^{\infty}(\G)$, the space of real valued and bounded functions on $\G$, {\it i.e.}, a linear form
%$
%m: \ell^{\infty}(\G) \longrightarrow \R
%$
%such that:
%\begin{itemize}
%\item[1.] $m( c) =c$, for a constant function $c$;
%\item[2.] $m(\f_1)\geq m(\f_2)$,  if $\f_1(\g)\geq \f_2(\g)$ for every $\g\in \G$;
%\item[3.] $m(\g_*\f)= m(\f)$, where for $\g\in \G$ and for $\f: \G\longrightarrow \R$, the function $\g_*\f$ is defined by
%$g_*\f(\g')=\f(\g'\g)$ for each $\g\in \G$.\\
%\end{itemize}
\end{Rem}

\noindent{\bf III.}
%Even though we shall not be interested in this case, 
For the sake of completeness let us mention that the polinomial growth condition is essentially optimal to have both existence and uniqueness of the asymptotic cone.
% For {\bf arbitrary finitely generated group}, L. van den Dries and A. J. Wilkie  \cite{vDW} generalised Gromov's construction: their sophisticated construction is based on the use of the theory of {\it non-standard extentions} and involves the choice of a  {\it non-principal ultrafilter}, which plays a crucial role in the determination of the asymptotic cone. In fact, S. Thomas and B. Velickovic exhibited in \cite{TV}  an example of a finitely generated group $\G$ and two non-principal ultrafilters $\omega$ and $\omega'$, such that the asymptotic cones $\cone^{\omega}(\G)$ and $\cone^{\omega'}(\G)$ are not homeomorphic to each other. 
We refer the interested readers to \cite{vDW, TV} and references therein.

\vspace{10 pt}

\subsection{Asymptotic Cone of Nilpotent Groups} \label{secasympnil}

In the case of a finitely generated (virtually) nilpotent group, it is possible to provide a very precise description of the asymptotic cone and of its properties.
In \cite{Pansu}, in fact, Pierre Pansu proved the following theorem (see also \cite{Gromov} and also the very nice presentation in \cite{BL}, where the authors, amongst other things, discuss the rate of this convergence):\\

\begin{Theorem}[Pansu] 
If $\G$ is a finitely generated virtually nilpotent group, then the sequence of metric spaces $\G_\e=(\G, d_\e)$ -- for any word or orbit metric -- converges in the pointed Gromov--Hausdorff topology 
 to a limit space  $(G_\infty, d_\infty)$, where $G_\infty$ is a connected, simply connected nilpotent (graded) Lie group and $d_\infty$ is a left-invariant Carnot-Caratheodory metric on $G_\infty$.\\
 \end{Theorem}

\begin{Rem}
The limit Lie group $G_\infty$ and the limit metric $d_{\infty}$ only depend on $\Gamma$.\\
\end{Rem}

In the following subsections, we aim to explain more carefully this result and the objects that are mentioned  in its statement.\\

%%%%%%%%%%%%%%%%%%%%%%%%%%%%%%%%%%%%%%%%%%%%%%%%%%%%%%%%%%
%%%%%%%%%%%%%%%%%%%%%%%%%%%%%%%%%%%%%%%%%%%%%%%%%%%%%%%%%%
%%%%%%%%%%%%%%%%%%%%%%%%%%%%%%%%%%%%%%%%%%%%%%%%%%%%%%%%%%

\subsubsection{Nilpotent Lie groups} \label{secstratified}
Let $\G$ be a finitely generated nilpotent group with no torsion elements. We recall that a group $\G$ is said to be {\it nilpotent} if its lower central series terminates in the trivial subgroup after finitely many steps:
\begin{equation}\label{centralseries}
\Gamma^{(1)}:=\G \geq \Gamma^{(2)}:=[\Gamma^{(1)}, \G] \geq \ldots \geq
\Gamma^{(i+1)}:=[\Gamma^{(i)}, \G] \geq \ldots \geq \G^{r} > \G^{r+1}=\{e\};
\end{equation}
the smallest $r$ such that $\G$ has a lower central series of length $r$ is called the nilpotency class of $\G$.\\

\noindent {\bf Examples.} \begin{itemize}
\item every abelian group is nilpotent (with $r=1$);
\item the smallest non-abelian example is provided by the {\it quaternion group} $Q_8$ (with $r=2$);
\item the direct product of nilpotent groups is itself nilpotent;
\item every finite nilpotent group is direct product of $p$-groups (which are nilpotent);
\item the Heisenberg group $\H_{2n+1}(\Z)$ is an example of infinite non-abelian nilpotent group (of nilpotency class $2$). Recall that this group is defined by
$$
H_{2n+1}(\Z) = \left\langle
a_1,b_1,\ldots, a_n,b_n,t\; \big|\; [a_i,b_i]=t\; \mbox{for each $i$ and all other brackets are $0$}
\right\rangle.
$$
%One can check that its nilpotency class is $2$.\\
\end{itemize}

\vspace{10 pt}

\begin{Rem}\label{canonicalgen}
Wolf in \cite{Wolf} proved that a finitely generated nilpotent group has polinomial growth, \ie there exist $C>0$ and $K>0$ such that
for each $r>0$
$$
\sharp \{\g\in \G: \; d(e,\g)\leq r \}  \leq C r^K
$$
(this notion is independent of  the metric $d$). It is easy to check that if $\Gamma$ is abelian and we choose the word metric on $\Gamma$, then $K$ coincides with the $\rank \G$. More generally, 
Bass proved in \cite{Bass} that if $\Gamma$ is a finitely generated nilpotent group, then
\begin{equation}\label{homdim}
K=
\sum_{k=1}^{r} k\cdot {\rm rank}\,\left( \Gamma^{(k)}/\Gamma^{(k+1)}\right),
\end{equation}
sometimes called the {\it homogeneous dimension} of the group.
%To understand where this unusual expression comes from, let us consider this example. Suppose $\G$ is generated by $\g_1,\ldots, \g_n \in \G$ and nilpotent of class 2 (this is the case of the above-mentioned Heisenberg group). Because
%commutators commute with all the generators, one may reduce words to a normal
%form. For instance, one can normalize $\g = \g_1 \g_2 \g_1 \g_3^{-1}$ as
%$$
%\g=\g_1 [\g_1 ,\g_2 ][\g_1 ,\g_2]^{-1}\g_2 \g_1 \g_3^{-1} =\g_1^2\g_2\g_3^{-1}[\g_1,\g_2]^{-1}.
%$$
%More generally, any word of length $r$ can be put in normal form in at most $r^2$ transpositions of adjacent elements. 
%Suppose the commutator subgroup $\G^{(1)}$ is generated by elements $h_1, \ldots, h_k$, then a word of length $r$ has a normalized form
%$\g = \g_1^{a_1} \cdot\ldots \cdot \g_m^{a_m} h_1^{b_1} \cdot \ldots \cdot h_k^{b_k}$
%where $a_i \leq r$. The normalization introduces at worst $r^2$ commutators, which  are then represented using $h_1, \ldots, h_k$, introducing a multiplicative term independent of $r$, so $\sum_i b_i \leq r^2$. Combining these estimates, we have 
%$$
%\sharp \{\g\in \G: \; d(e,\g)\leq r \}  \leq \tilde{C} r^{n+2k}.\\
%$$
\end{Rem}

%\vspace{10 pt}

\subsubsection{Malcev's Closure} The importance of nilpotent groups, as far as this discussion is concerned, comes from the following result. In \cite{Malcev}  Malcev proved that every finitely generated torsion-free nilpotent group $\Gamma$ can be ``completed" to be a simply connected nilpotent Lie group, which is often called the {\it Malcev's closure} of $\Gamma$ (exactly in the same way $\R^n$ can be thought as a completion of $\Z^n$). More precisely,

\begin{Theorem}[Malcev] Let $\G$ be a finitely generated torsion free nilpotent group. There exists a unique (up to isomorphisms) simply connected nilpotent Lie group $G$, in which $\G$ can be embedded as a co-compact discrete subgroup.
\end{Theorem}

%The idea behind Malcev's construction is essentially the following. First, Malcev proved that every finitely generated nilpotent torsion-free group $\Gamma$ possesses a canonical set of generators; namely, there exist elements $\g_1, \ldots, \g_d \in \G$ such that every $\g\in\G$ may be written uniquely as $\g = \g_1^{a_1}\cdot \ldots \cdot \g_d^{a_d}$ with $a_1,\ldots, a_d\in \Z$ (the idea is similar to the one explained in Remark \ref{canonicalgen}). Then, he proves that the group structure on $\Gamma$ can be defined in terms of polynomials in the $a_i$'s; these polynomials - if one considers $a_i\in\R$ rather than integers - define a group structure on $\R^d$ which defines the completion group $G$.\\

\begin{Rem}\label{rem10}\noindent
\begin{itemize}
\item The dimension of $G$ is given by  
$
\sum_{k=1}^{r} {\rm rank}\,\left( \Gamma^{(k)}/\Gamma^{(k+1)}\right).
$
\item If $\G=\Z^n$ this construction gives exactly $G=\R^n$. As a less-trivial example, one can consider the Heisenberg group
$$
H_{2n+1}(\Z) = \left\langle
a_1,b_1,\ldots, a_n,b_n,t\; \big|\; [a_i,b_i]=t\; \mbox{for each $i$ and all other brackets are $0$}
\right\rangle.
$$
It is not difficult to check that in this case Malcev's construction leads to
$
G=H_{2n+1}(\R) \simeq \C^n\times \R$,
with the group structure given by $(z,t)*(z',t')=(z+z',t+t'+{\rm Im}\langle z,z'\rangle)$,
where $\langle \cdot,\cdot \rangle$ denotes the standard Hermitian pairing of complex vectors.\\
\end{itemize}
\end{Rem}

As we are going to describe in the following subsection, an important property of this completion is that to every simply connected nilpotent Lie group $G$ one can associate a {\it graded Lie group} (or {\it Carnot group}) $G_\infty$. In particular, this new group will admit a one-parameter subgroup of $\R$-diagonalisable automorphisms that can be interpreted as dilations (or homotheties), and that will be of fundamental importance in the homogenization process. See also \cite{LeDonne} for a nice presentation of these objects. \\

\subsubsection{Stratified Lie Algebras and Graded Lie Groups} \label{subsubseccarnotgroup}
A Lie algebra $\fg$ is called  a {\it stratified algebra} if it admits a stratification, namely there exist vector subspaces $V_1,\ldots, V_r \subset \fg$ (called {\it strata}) such that
$$
\fg = V_1 \oplus \ldots \oplus V_r
$$
and $[V_j,V_1]=V_{j+1}$ for $j=1,\ldots, r-1$, with $V_r\neq \{0\}$ but $[V_r,V_1]=\{0\}$ (by $[V,W]$ we mean the vector subspace generated by commutators of the form $[v,w]$ with $v\in V$ and $w\in W$). \\

A  Lie group is said to be {\it graded} if it is simply connected and its Lie algebra is stratified. \\

\begin{Rem}\noindent
\begin{itemize}
\item A stratified Lie algebra is nilpotent with nilpotency class equals to the number of strata. 
\item The first stratum $V_1$ completely determines the other strata.
\item The commutator subalgebra $[\fg,\fg]$ coincides with $V_2\oplus \ldots \oplus V_r$ and therefore $V_1$ is in bijection with the the abelianization $\frac{\fg}{[\fg,\fg]}$.\\
\end{itemize}
\end{Rem}

Let us go back to the core of our discussion. Although the Malcev's closure $G$ of a finitely generated nilpotent group $\G$ is not necessarily  graded, yet there is a canonical way to associate to it a graded algebra.
Let $\fg$ be a nilpotent Lie algebra of nilpotency class $r$ and let us consider the descending central series given by $\fg^{(1)}:=\fg$ and, inductively, $\fg^{(i+1)}:=[\fg, \fg^{(i)}]$ for $i=1,\ldots, r$. The {\it graded Lie algebra} associated to $\fg$ is given by the Lie algebra $\fg_\infty$ given by the direct-sum decomposition
$$
\fg_{\infty} := \bigoplus_{i=1}^r \fg^{(i)}/\fg^{(i+1)}
$$
endowed with the unique Lie bracket $[\cdot,\cdot]_\infty$ that has the property that if $v\in \fg^{(i)}$ and $w\in \fg^{(j)}$, then the bracket is defined as 
$$
[\overline{v},\overline{w}]_\infty = \overline{[v,w]} \qquad  ({\rm modulo}\; \fg^{(i+j+1)}).
$$

Notice that $\fg_{\infty}$ is a stratified algebra. One can then consider the unique connected, simply connected Lie group $G_\infty$ whose Lie algebra is $\fg_\infty$. This group will be called the {\it graded group} of $\fg$.\\

\begin{Rem}\label{algebras}
One can identify $\fg$ and $\fg_{\infty}$ (yet, this identification is not unique, as it will be clear from the following). Let us choose $r$ subspaces of $\fg$, namely 
$V_1,\ldots, V_r$ such that $\fg^{(j)}= \fg^{(j+1)}\oplus V_j$ for each $j=1,\ldots, r$. In particular, the projection whose kernel is $g^{(j+1)}$ provides a linear isomorphism between
$V_j$ and $\fg^{(j)}/\fg^{(j+1)}$ and this induces a linear isomorphism between $\fg$ and $\fg_{\infty}$. In particular, identifying $G$ with $\fg$ and $G_\infty$ with $\fg_\infty$
via the respective exponential maps, one finds an identification between $G$ and $G_\infty$; this identification is a diffeomorphism, but not a group isomorphism. However, one can use this identification to pull back  the lie group structure of $G_\infty$ on $G$, and define a new Lie product on $G$ which makes it a graded Lie group.\\
\end{Rem}

\subsubsection{Metric Structures on Graded Lie Groups: Carnot Groups}\label{subsubseccarnotgroup2}
Let $G$ be a graded Lie group with 
stratified algebra
$$
\fg = V_1 \oplus \ldots \oplus V_r.
$$
Let us consider a norm $\|\cdot\|$ on $V_1$ and let  $\Delta$  denote the horizontal distribution induced by $V_1$, namely the left-invariant subbundle of the tangent bundle of $G$ such that a vector $v$ at a point $p\in G$ is an element of $\Delta$ if and only if $L_p^*v \in V_1$ (where $L_p$ denotes the left multiplication by $p$). This allows one to extend the norm on $V_1$ to the whole $\Delta$, simply by defining $\|v\|:=\|L_p^*v\|$.\\

The triple $(G,\Delta, \|\cdot\|)$ is an example of subFinsler manifold. We can define on it a distance function, called {\it Carnot-Carath\'eodory-Finsler metric},  in the following way. A curve $\gamma:[0,t]\longrightarrow G$ is said to be {\it horizontal} (with respect to $\Delta$) if it is absolutely continuous and its velocity $\dot{\g}(s) \in V_1$ for almost all $s\in [0,t]$.
In particular, if $p,q \in G$ we define:
$$
d_{\rm sF}(p,q) := \inf \left\{
\int_0^t \|\dot{\g}(s)\|\,ds\,:\; \g\; \mbox{is horizontal from $p$ to $q$ and  $t>0$}
\right\}.
$$

Since, as remarked above,  the first stratum $V_1$ generates the whole $\fg$, then this function $d_{\rm sF}$ is finite (\ie every two points in $G$ can be connected by horizontal paths) and defines a geodesic distance that induces the topology manifold on $G$ (this follows from a theorem by Chow, see for example \cite[Chapter 2]{Mont}).\\

The group $G$ with the metric $d_{\rm sF}$ is called a {\it Carnot group}. A peculiarity of Carnot groups is that they admits dilations. In fact, for each $\l\in \R$, one can introduce the algebra-dilations $\delta_\l: \fg \longrightarrow \fg$ which are defined linearly by imposing $\d_\l(v)=\l^i v$ for every $v\in V_i$ with $i=1,\ldots, r$; if $\l\neq 0$ these are algebra automorphisms and 
the maps $\l \longmapsto \delta_\l$ for $\l\geq 0$ constitute a one-parameter subgroup of ${\mathcal A}{\rm ut}(\cg)$.\\
Observe that since $G$ is simply connected, than the exponential map is a diffeomorphism and therefore these dilations
induce unique automorphisms of $G$ (that we shall refer to as {\it group dilations} and continue to denote by $\d_\l$), which are given by
$\d_\l (p) = \exp \circ \delta_\l \circ \exp^{-1} (p)$ for all $\l\in \R$ and $p\in G$.
In particular, one can easily check that
$$
d_{\rm sF}(\delta_\l(p), \d_\l(q)) = \l\, d_{\rm sF}(p,q) \qquad \forall\;p,q\in G\quad {\rm and} \quad \l\in\R_+.\\
$$

In our case, on $G_\infty$ - \ie the graded Lie group associated to $\Gamma$ - we consider  the Carnot-Caratheodory metric $d_\infty$ associated to:
\begin{itemize}
\item a horizontal subspace $V_1$, which is transverse to the commutator subalgebra (see subsection \ref{subsubseccarnotgroup});
\item a norm $\|\cdot\|_\infty$ on $V_1$ defined in the following way (see \cite[Section 21]{Pansu} for more details). Let $|\cdot|$ denote a norm on $\Gamma$ (see subsection \ref{metrics}); this induces a norm on the abelian subgroup $A:=\frac{\G}{[\G,\G]}$, which can be seen a subspace of the vector space $B:=\frac{\cg}{[\cg,\cg]} =\frac{\cg^{(1)}}{\cg^{(2)}}$. If one homogenizes the norm on $A$, \ie considers\footnote{Compare this with what discussed in Case I in section \ref{asymptoticcone}. Moreover, this construction is closely related to the so-called {\it stable norm} on $H_1(M;\R)$, where $M$ is a Riemannian manifold,  defined by Federer and Fleming in \cite{FF}.} 
$$
|a|_{\infty} = \lim_{k\rightarrow +\infty} \frac{|k\cdot a|}{k},
$$
then this function can be extended to norm on $B$, making it a Banach space. Observe now that 
the homomorphism 
$G_{\infty}\longrightarrow \frac{\cg_{\infty}}{[\cg_{\infty}, \cg_{\infty}]}$ induces a isomorphism between $V_1$ and $B$, that can be used to transport the norm on $V_1$.\\
%
%
%
%\stackrel{\simeq}{\longrightarrow} V_1$. As proved by Malcev in \cite{Malcev}, the commutator subgroup $\G^{(2)}=[\G,\G] =\G \cap [G,G]$ and the image $\Pi (\Gamma)$ is a co-compact lattice in 
%$V_1$. The norm $\|\cdot\|_\infty$ on $V_1$ corresponds to the norm whose unit ball coincides with the convex envelope of elements of the form $\frac{\Pi(\g)}{|\g|}$, with $\g\in \G$.
\end{itemize}

\section{Homogenization Procedure and Proof of the Main Result} \label{Hom} 
In this section we shall prove the main result stated in the Introduction, namely a Homogenization theorem/procedure for Hamilton-Jacobi equation associated to a Tonelli Hamiltonian which is invariant under the action of a discrete nilpotent group.

\subsection{Notation and Assumptions} \label{notation}  Before proceeding with our discussion, in order to help the reader, we summarize and recall our notation and assumptions so far, and add a few more (simplifying) ones.\\

\begin{itemize}
\item  $X$ is a smooth connected (non-compact) manifold without boundary, endowed with a complete Riemannian metric $d$. $TX$ and $T^*X$ denote the tangent and cotangent bundles. We shall denote by $X_\e$ the rescaled metric space $(X,d_\e:=\e d)$. We also fix a reference point $x_0\in X$ (our main result will be independent of this choice).\\

\item $\G$ is a finitely generated torsion-free nilpotent group of nilpotency step $r$ and 
$|\cdot|$ denote a norm on $\Gamma$: either the algebraic  one $\|\cdot\|_S$ associated to a symmetric generating set $S$ or the one associated to the orbit metric $|\gamma|:=d_{\G,x_0}(x_0,\gamma\cdot x_0)$; we denote by $d_\Gamma$ the induced metric  (see subsection \ref{metrics}). $B_\G(r)$ will denote the closed ball of radius $r$ in $(\Gamma, |\cdot|)$ centered at the identity $e$.
\\

\item We assume that $\Gamma$ acts on $X$ by isometries (\ie $d$ is invariant under the action of $\Gamma$) and that this action is free, properly discontinuous and co-compact. As we have recalled in subsection \ref{sec1.2}, this action  naturally extends to actions on   $TX$  and $T^*X$.  \\

\item We denote by $G$ the Malcev's closure of $\Gamma$ and by $\cg$ its Lie algebra, which is nilpotent of class $r$; its descending central series will be denoted by $\fg^{(1)}:=\fg$ and, inductively, $\fg^{(i+1)}:=[\fg, \fg^{(i)}]$ for $i=1,\ldots, r$.\\

\item  $G_{\infty}$ represents the graded group associated to $G$, as explained in section \ref{secasympnil}, and $\cg_{\infty}$ is its stratified Lie algebra (with strata $V_1, \ldots, V_r$).
We shall denote with a bar ({\it e.g.},  $\bar{x}$) the elements of $G_\infty$.\\
Recall that  $G_\infty$ comes equipped with a Carnot-Carath\'eodory distance $d_\infty$ and a one-dimensional family of dilations $\delta_t$ (which can be interpreted, via the exponential map, either as automorphisms of the group or of the algebra). We shall indicate with $B_\infty(R)$ the closed ball of radius $R$ in $(G,d_\infty)$ centered at the identity $e$. 
\\

\item  We fix an identification between $\cg$ and $\cg_{\infty}$ (see also Remark \ref{algebras}).
Let us choose $W_i$ a supplementary subspace of $\cG^{(i+1)}$ in $\cG^{(i)}$, \ie $\fg^{(i)}= \fg^{(i+1)}\oplus W_i$  and  consider the unique linear bijection $L: \cG \longrightarrow \cG_{\infty}$ such that 
%\begin{itemize}
%\item[a)] 
$L(W_i)=V_i$, where $V^i$ is the $i$-th stratum in the stratification of $\cG_{\infty}$,
and -- modulo the isomorphisms $W_i\simeq \cG^{(i)} / \cG^{(i+1)}$ and $V_i \simeq \cG^{(i)} / \cG^{(i+1)}$ --  $L$ corresponds to the identity map.\\
%\end{itemize}

\item $H: T^*X \longrightarrow \R$ is a $\G$-{invariant  Tonelli} Hamiltonian and  $L: T X \longrightarrow \R$ is the corresponding $\G$-{invariant  Tonelli} Lagrangian.\\

\item
Up to modifying the Lagrangian by adding a constant, we can assume without loss of generality that $\beta(0)=-\min_{c\in H^1( X\slash\Gamma;\R)} \a(c) =0$, where $\a=\a_L$ and $\b=\b_L$ denote Mather's minimal average actions associated to $L:  T( X\slash\Gamma) \longrightarrow \R$ (hereafter we omit the dependence on $L$ to simplify the notation). See Appendix \ref{appMather} and  \cite{Mather91, Gonzalobook, SorrentinoLectureNotes}.\\

\item In order to simplify the presentation of the proof, we also assume that $H_1(X;\R)=0$, so that $\frac{\cg}{[\cg,\cg]} \simeq H_1( X\slash\Gamma;\R)$. Otherwise, all results  remain true, up to identifying $\frac{\cg}{[\cg,\cg]}$ with a (possibly proper) subgroup of $H_1( X\slash\Gamma;\R)$.\\

\end{itemize}

%\begin{Rem}\label{remvirtnilp}
%One could also consider the more general case of $\G$ being a finitely generated torsion-free group with polinomial growth; then, by the already-mentioned theorem by Gromov \cite{Gromov} (see also subsection \ref{asymptoticcone}), $\G$ is virtually nilpotent, \ie it contains a nilpotent torsion-free subgroup $\Gamma'$ of finite index.  The strategy would be to consider the action of $\Gamma'$. In fact, the torsion elements  form a finite subgroup $T$ in $\Gamma$ (see \cite[Proposition 67]{Pansu}) and one can consider the quotient group $\Gamma'=\Gamma/T$, which is nilpotent and torsion-free; observe that since $\Gamma$ and $\Gamma'$ have finite Hausdorff distance from each other, then they have the same asymptotic cone.\\ \end{Rem}

\subsection{Outline of the Proof}
\begin{itemize}
\item[-] In subsection \ref{secrescHJe} we consider rescaled Hamilton-Jacobi equations and discuss properties of their solutions.
\item[-] In subsection \ref{secrescmap} we introduce a useful notion of convergence for functions defined on rescaled metric spaces:  this  will be  achieved by defining suitable rescaling maps that will allow us to ``transport'' these functions to/from the limit space.
\item[-] In subsection \ref{secgenbeta} we study the convergence of the Ma\~n\'e potentials associated to the rescaled Lagrangians and define what we shall call ``generalized Mather's $\beta$-function'' (in accordance to what happens in the abelian case). This map plays the role of the effective Lagrangian for the limit problem.
\item[-] Finally, in subsection \ref{secproofmain}  we prove the convergence of solutions to the rescaled problem to a solution to the limit problem: this will complete the proof of our main result.
\end{itemize}
\vspace{10 pt}

%%%%%%%%%%%%%%%%%%%%%%%%%%%%%%%%%
%\subsection{{\bf Rescaled Hamilton-Jacobi equation}}
%
%We want now to consider the following problem ($\e>0$):
%
%
%\begin{equation} \label{HJeps}
%\left\{ \begin{array}{l}
%\partial_t v^\e(x,t) + H(x, \frac{1}{\e} \partial_x v^\e(x,t)) = 0 \qquad x\in X,\; t>0\\
%v^\e(x,0)= f_\e( x),\\
%\end{array}\right. 
%\end{equation}
%
%
%\noindent where $f_\e: X_\e \longrightarrow \R$ are uniformly continuous functions that converge locally uniformly to  a function $f: G_\infty \longrightarrow \R$ which has at most linear growth \textcolor{red}{{\bf (IPOTESI??? $f_\e$ locally equiLipschitz?)}}. Note that this problem corresponds to Hamilton-Jacobi equation with a rescaled metric $d_\e$ on $X$.
%Let us denote by 
%$ H_\e(x,p) := H(x, \frac{1}{\e} p)$
%the  new Hamiltonian and by  $L_\e(x,v) := L(x, {\e} v)$ the corresponding Lagrangian.

%%%%%%%%%%%%%%%%%%%%%%%%%%%%%%%%%

\subsection{{Rescaled Hamilton-Jacobi Equation}}\label{secrescHJe}

Given $\e>0$, we consider the following problem:

\begin{equation} \label{HJeps}
\left\{ \begin{array}{l}
\partial_t v^\e(x,t) + H(x, \frac{1}{\e} \partial_x v^\e(x,t)) = 0 \qquad x\in X_\e,\; t>0\\
v^\e(x,0)= f_\e( x),\\
\end{array}\right. 
\end{equation}

\noindent where $f_\e: X_\e \longrightarrow \R$ are equiLipschitz functions on  $X_\e=(X, d_\e:=\e d)$. As we have explained in subsection \ref{homgenmanif},  this problem corresponds to Hamilton-Jacobi equation with a rescaled metric $d_\e$ on $X$.
We would like to give a meaning to the fact  that these functions  converge to some  function  $\bar{f}$ defined on a {\it limit space} ($G_\infty$ would be the best candidate); we shall discuss this issue in subsection \ref{subsecheps}.\\

Let us now start by studying the dependence of this problem and of its solutions on $\e$.
Let us denote by 
$ H_\e(x,p) := H(x, \frac{1}{\e} p)$
the  new Hamiltonian and by  $L_\e(x,v) := L(x, {\e} v)$ the corresponding Lagrangian.\\

It is a classical result that the variational solution to (\ref{HJeps}) is given by the Lax--Oleinik  formula (see for example \cite{Evans, Fathibook, Lions} ):

\begin{eqnarray}\label{vepsilon}
v^\e(x,T) &=& \inf  \left\{
v^\e(\g(0),0) + \int_0^T L_\e (\g(t),\dot{\g}(t))\,dt\; \big| \; \g\in C^1([0,T], X),\; \g(T)=x \right\}\nonumber\\
&=&
\inf  \left\{
v^\e(\g(0),0) + \int_0^T L (\g(t), \e \dot{\g}(t))\,dt\; \big| \; \g\in C^1([0,T], X),\; \g(T)=x \right\}.
\end{eqnarray}

If we consider the curve $\eta:\left[0,\frac{T}{\e}\right] \longrightarrow X$, given by $\eta(s):=\gamma(\e s)$, then
$$
\int_0^T L(\gamma(t), \e \dot{\g}(t))\, dt = \e \int_0^{\frac{T}{\e}} L(\eta(s), \dot{\eta}(s))\,ds,
$$
and therefore:
\begin{eqnarray*}
v^\e(x,T) &=& 
\inf  \left\{
v^\e(\g(0),0) + \e \int_0^{{T}/{\e}} L(\eta(s), \dot{\eta}(s))\,ds,
\; \big| \; \eta\in C^1([0,{T}/{\e}], X),\; \eta({T}/{\e})=x \right\} \\
&=&
\inf_{y\in X} \left\{ v^\e(y,0) + \e\, \phi(y,x,T/{\e}) \right\}, 
\end{eqnarray*}

where $\phi$ represents  the so-called Ma\~n\'e potential (see for instance \cite{Gonzalobook, SorrentinoLectureNotes} for more details):
\begin{equation}\label{manepotential}
\phi(y,x,S) := \inf \left\{ \int_0^S L (\eta(s),\dot{\eta}(s))\,ds\; \big| \; \eta\in C^1([0,S], X),\; \eta(0)=y, \;\eta(S)=x \right\}\!;
\end{equation}
this quantity is well-defined since Tonelli's theorem  (on the existence of action-minimizing curves) holds in the non-compact setting as well, as long as the Lagrangian is superlinear  (see \cite[Chapter 3]{Gonzalobook}). 

Moreover, one can prove the following result.\\ %(we adapt here the proof from \cite[Lemma 3.1]{CIS}).\\

\begin{Lem}\label{lemma1}
Let $\lambda>0$. Then, there exists a positive constant $K=K(\lambda)>0$ such that for each $T\geq \lambda$ and for each $x,y,z\in X$
$$
d(y,z)<\frac{T}{2} \qquad \Longrightarrow \qquad
|\phi(x,y,T)- \phi(x,z,T)| \leq K d(y,z).\\
$$
\end{Lem}

%\begin{Lem}\label{lemma1}
%Let $\lambda>0$. Then, there exists a positive constant $K=K(\lambda)>0$ such that for each $T\geq \lambda$  the map $(x,y) \longmapsto \phi(x,y,T)$ is Lipschitz on the set $\{(x,y): \; d(x,y)<T\}$.\\
%\end{Lem}

\begin{Rem}
The Lipschitz constant $K=K(\l)$ that one obtaines from the proof goes to infinity as $\lambda$ goes to zero. \\
\end{Rem}

\begin{Proof}
Let $T\geq \lambda$. It follows from the superlinearity of the Lagrangian that there exists $A=A(\lambda)>1$ such that if $\g: [0,T] \longrightarrow X$ is a Tonelli minimizer, %and $d(\g(0), \g(T))<T$, 
then
$\|\dot{\g}(t)\|\leq A$ for all $t\in [0,T]$ (see \cite[Appendix 1]{Mather91} and \cite[Lemma 3.2-1 \& Corollary 3.2-2]{Gonzalobook}).
Let 
$$
Q_1:= \sup_{\|v\|\leq 2A} |L(x,v)| \qquad {\rm and} \qquad
Q_2:= \sup_{\|v\|\leq 2A} |\partial_v L(x,v) \cdot v|.
$$
 If $a\in[\frac{1}{2},2]$ and $\|v\|\leq A$ we have that:
\begin{eqnarray}\label{stimettina}
\left|
L(x,av)\cdot \frac{1}{a} - L(x,v) \right| &\leq&
\frac{1}{a} \left| L(x,av) - L(x,v) \right| + \left|\frac{a-1}{a}\right| |L(x,v)| \nonumber\\
&\leq&
\frac{1}{a} \left( \int_1^a \partial_v L(x,sv)\cdot v\, ds \right) + \left|\frac{a-1}{a}\right| |L(x,v)| \nonumber \\
&\leq&
\left|\frac{a-1}{a}\right| Q_2 + \left|\frac{a-1}{a}\right| Q_1.
\end{eqnarray}
Let now $z\in X$ such that $d:=d(y,z)< \frac{T}{2}$ and let $\g:[0,T] \rightarrow X$ be a Tonelli minimizer such that $\g(0)=x$ and $\g(T)=z$, {\it i.e.},
$$
\phi(x,z, T) = \int_0^T L(\g(t), \dot{\g}(t))\,dt.
$$
Define a new curve $\eta: [0, T-d] \longrightarrow X$ by $\eta(s):= \g \left( s\cdot \frac{T}{T-d}\right)$. Then, applying (\ref{stimettina}) with $a= \frac{T}{T-d}$ (which, because of our choice of $T$, belongs to $[1/2,2]$):
\begin{eqnarray*}
\phi(x,z, T-d) &\leq& 
\int_0^{T-d} L(\eta(s), \dot{\eta}(s))\, ds =
\int_0^{T} L\left(\g(t), \frac{T}{T-d}\dot{\g}(t)\right) \cdot \left(\frac{T-d}{T} \right)\, dt \\
&\leq&
\int_0^{T} L\left(\g(t), \dot{\g}(t)\right)\, dt  + \frac{d}{T-d}  (Q_1+Q_2) \\
&\leq& 
\phi(x,z,T) + \frac{2}{T} (Q_1+Q_2) \, d(y,z). \\
\end{eqnarray*}
Moreover,  considering the shortest geodesic connecting $z$ to $y$, one also obtains:
$$
\phi(z,y,d) \leq Q_1 d(y,z).
$$
Using the triangle inequality, we conclude that:
\begin{eqnarray*}
\phi(x,y,T) &\leq& \phi(x,z,T-d) + \phi(z,y,d) \leq
\phi(x,z,T) + \left[ Q_1 + \frac{2}{T}(Q_1+Q_2)\right]\,d(y,z).
\end{eqnarray*}
Similarly, one proves that
$$
\phi(x,z,T) \leq \phi(x,y,T) + \left[ Q_1 + \frac{2}{T}(Q_1+Q_2)\right]\,d(y,z).
$$
This concludes the proof of the theorem with $K:=Q_1 + \frac{2}{\lambda}(Q_1+Q_2)$.\\
\end{Proof}

Using this lemma, one can now prove the following result.

\begin{Prop}\label{equilip}
Let $r>0$. There exists $K=K(r)>0$ such that if $T\geq 2 \,r \e$, then $v^{\e}(x,T)$ is $\e K$-Lipschitz (in the $x$-component) in any ball of radius $r$  (with respect to the metric $d$).\\
% (more precisely, on balls of radius $C$ with respect to the metric $d$)
Moreover, if $\e<1$ and $T\geq 2\,r$ then the function $v^{\e}(\cdot, T): X_\e \longrightarrow \R$ is  $K$-Lipschitz (in the $x$-component) in any ball of radius $r$ (with respect to the rescaled metric $d_\e$).\\
\end{Prop}

\begin{Rem}
The Lipschizt constant $K=K(r)$ can be chosen to  be the same as the one in Lemma \ref{lemma1} (with $\lambda=r$).\\
\end{Rem}

%\begin{Prop}\label{equilip}
%Let $C>0$. Then there exists $K_1=K_1(C)>0$ such that if $\frac{T}{\e}>C$, then $v^{\e}(x,T)$ is $\e K_1$-Lipschitz in the $x$ component  in any ball of radius $C$ (with respect to the metric $d$).
%% (more precisely, on balls of radius $C$ with respect to the metric $d$)
%In particular, if $\frac{T}{\e}>C$ the function $v^{\e}(\cdot, T): X_\e \longrightarrow \R$ is locally $K_1$-Lipschitz in any ball of radius $C$ (with respect to the rescaled metric $d_\e$).\\
%\end{Prop}
%

\begin{Proof}
%The proof is the same as \cite[Lemma 3.2]{CIS}, we rewrite it here for the reader's convenience.
Let $x,z\in X$ with $d(x,z)<r<\frac{T}{2\e}$. Let $K=K(r)$ be the Lipschitz constant from Lemma \ref{lemma1}
and let $y_n\in X$ be such that
$$
v^{\e}(x,T)= \lim_{n\rightarrow +\infty} \left\{ v^\e(y_n,0) + \e \phi(y_n,x, {T}/{\e}) \right\}.
$$
Then, observing that  $\phi(y_n,x, {T}/{\e}) \geq \phi(y_n,z, {T}/{\e}) - K d(x,z)$, 
we conclude that:
\begin{eqnarray*}
v^{\e}(x,T) &=& \lim_{n\rightarrow +\infty} \left\{ v^\e(y_n,0) + \e \phi(y_n,x, {T}/{\e}) \right\} \\
&\geq& \limsup_{n\rightarrow +\infty} \left\{ v^\e(y_n,0) + \e \phi(y_n,z, {T}/{\e}) - \e Kd(x,z) \right\} \\
%&=& \limsup_{n\rightarrow +\infty} \left\{ v^\e(y_n,0) + \e \phi(y_n,z, {T}/{\e}) \right\} - \e Kd(x,z)  \\
&\geq& v^{\e}(z,T)  - \e K d(x,z),
\end{eqnarray*}
or equivalently
$$
v^{\e}(z,T)  - v^{\e}(x,T)  \leq \e K d(x,z).
$$
\noindent The reversed inequality can be proven similarly.\\

Let us now prove the second part of the statement. Let $x,z\in X_\e$ with $d_\e(x,z)<r$ and $T\geq 2\,r$ (recall that  $d_\e(x,z) := \e d(x,y)$). Therefore, $d(x,z) = \frac{d_\e(x,z)}{\e} < \frac{r}{\e} \leq \frac{T}{2\e}$.
 Let $K=K(r)$ be the Lipschitz constant from Lemma \ref{lemma1} (observe that since $\e<1$, then $r\leq \frac{r}{\e}$)
and let $y_n\in X_\e$ be such that
$$
v^{\e}(x,T)= \lim_{n\rightarrow +\infty} \left\{ v^\e(y_n,0) + \e \phi(y_n,x, {T}/{\e}) \right\}.
$$
Then, observing that  $\phi(y_n,x, \frac{T}{\e}) \geq \phi(y_n,z, {T}/{\e}) - K d(x,z)$, 
we conclude that:
\begin{eqnarray*}
v^{\e}(x,T) &=& \lim_{n\rightarrow +\infty} \left\{ v^\e(y_n,0) + \e \phi(y_n,x, {T}/{\e}) \right\} \\
&\geq& \limsup_{n\rightarrow +\infty} \left\{ v^\e(y_n,0) + \e \phi(y_n,z, {T}/{\e}) - \e Kd(x,z) \right\} \\
%&=& \limsup_{n\rightarrow +\infty} \left\{ v^\e(y_n,0) + \e \phi(y_n,z, {T}/{\e}) \right\} - \e Kd(x,z)  \\
&\geq& v^{\e}(z,T)  - \e K d(x,z) = v^{\e}(z,T)  - K d_\e(x,z),
\end{eqnarray*}
or equivalently
$$
v^{\e}(z,T)  - v^{\e}(x,T)  \leq K d_\e(x,z).
$$
\noindent The reversed inequality can be proven similarly.\\
\end{Proof}

%%%%%%%%%%%%%%%%

%\vspace{20 pt}

\subsection{Limit problem, Rescaling Maps and Convergence}\label{secrescmap}

In this subsection we want to discuss  the convergence to the limit problem, introducing a suitable notion of convergence.
For, we need to introduce suitable ``{\it rescaling maps}'' that will allow us to ``transfer'' solutions to the rescaled problems on the limit metric space (our inspiration came from \cite{Pansu}).

\subsubsection{Rescaling Maps} \label{subsecheps}
We distinguish two different cases (the first can be considered as a subcase of the second, but we believe it might be useful itself for illustrative purposes):
\begin{itemize}
\item[I.] $\Gamma$ is a finitely generated torsion-free nilpotent group which is already a discrete co-compact subgroup of a graded nilpotent Lie group $G=G_\infty$ ({\it e.g.}, abelian groups or the Heisenberg group, see subsection \ref{secstratified}).
\item[II.] $\Gamma$  is a general finitely generated torsion-free nilpotent group.\\
\end{itemize}

\begin{itemize}
\item[{Case I.}]
{\it The limit space.} Let us start with the easier case in which $\Gamma$ is a finitely generated nilpotent group which is already a discrete co-compact subgroup of a graded nilpotent Lie group $G=G_\infty$, with stratified Lie algebra $\cG=\cG_{\infty}$. \\
%As we have discussed in Section \ref{secasympnil} (see also \cite{Pansu}), $G$ is equipped with a Carnot-Carath\'eodory distance $d_\infty$ and a one-dimensional family of dilations $\delta_t$, which are group automomorphisms; these dilations can be also lifted to algebra homomorphisms of $\cG$ via the inverse of the exponential map (which is a global diffeomorphism, since $G$ is simply connected); we shall continue to denote them by $\delta_t$ (it will be clear from the context which we are considering).

%
%Let us fix $x_0\in X$ and let $|\cdot|$ denote a norm on $\Gamma$: either the algebraic  one $\|\cdot\|_S$ associated to a symmetric generating set $S$ or the one associated to the orbit metric $|\gamma|:=d_{\G,x_0}(x_0,\gamma\cdot x_0)$ (see subsection \ref{metrics}).\\

\noindent{\it Construction of the rescaling maps.} % \label{subsecheps}
Let us choose a compact fundamental domain $\Omega$ for the action of $\Gamma$ in $G$ (we can assume that the identity $e\in \Omega$);  for each $\e>0$ let us define the {\it rescaling maps} $h_\e : G \longrightarrow \Gamma$ such that for each $\bx\in G$, $h_\e(\bx)$ equals an element $\g\in \G$ such that $\delta_{1/\e}(\bx) \in \gamma \cdot \Omega$.\\

\begin{Rem}\label{propertyfecase1}
Properties of these maps $h_\e$  (see \cite[\S 27]{Pansu}):
\begin{itemize}
\item[i)]  For each $R>0$, there exists $\theta(\e)\rightarrow 0$ as $\e\rightarrow 0$ such that $h_\e$ maps the ball $B_\infty(R)$ into $B_\G\left(\frac{R+\theta(\e)}{\e}\right)$, which corresponds to the ball of radius $R+\theta(\e)$ with respect to the rescaled norm $\e|\cdot|$.
\item[ii)] If $\g_\e$ is a sequence in $\Gamma$ such that $\delta_\e \gamma_\e \rightarrow \bx \neq e$ as $\e\rightarrow 0$, then $\e |\gamma_\e|$ tends to $d_{\infty}(e,\bx)$ (see \cite[Proposition 41]{Pansu}).  In particular, for each $\bx\in G\setminus\{e\}$ we have that $\delta_\e(h_\e(\bx))$ converges to $\bx$ in G; in fact:
$$
d_\infty(\delta_\e(h_\e(\bx)), \bx) = \e d_\infty(h_\e(\bx), \delta_{1/\e}(\bx)) \leq C \e \,{\rm diam}(\Omega) \stackrel{\e\rightarrow 0^+}{\longrightarrow} 0.
$$
\item[iii)] Moreover,
if $\bx \neq \by$, then $\delta_\e\left( h_\e(\bx)^{-1} h_\e(\by)\right)$ tends to $\bx^{-1}\by$ in $G$ and, in particular, 
$\e |h_\e(\bx)^{-1} h_\e(\by)| \rightarrow d_\infty(\bx,\by)$ as $\e$ tends to zero; this convergence is uniform on compact sets of $G$.
\\
\end{itemize}
\end{Rem}

\item[{ Case II.}] {\it The limit space.}  Let us consider now the general case in which $\Gamma$ is a general finitely generated torsion-free nilpotent group. 
As we have recalled in section \ref{subsubseccarnotgroup}, $\Gamma$ embeds as a co-compact lattice in its Malcev closure $G$; however, differently from case I, this group  might not be {\it graded}, hence one needs to consider the associated  {\it graded group} $G_\infty$.
Recall that this group can  be naturally identified with $G$, but this identification is not a homomorphism (it is only a diffeomorphism); however, one could interpret this graded group as $G$ with a different Lie structure and a different group operation.\\

%$G_\infty$ is equipped with a Carnot-Carath\'eodory distance $d_\infty$ and a one-dimensional family of dilations $\delta_t$, which are group automomorphisms; these dilations can be also lifted to algebra homomorphisms of $\cG_\infty$ via the inverse of the exponential map $\exp_{\infty}$ (which is a global diffeomorphism, since $G$ is simply connected); we shall continue to denote them by $\delta_t$ (it will be clear from the context which we are considering).\\

\noindent{\it Construction of the rescaling maps.} \label{subsecheps2}
The main problem in constructing the rescaling maps in this case is that,  differently from case I,  on $G$  (or equivalently on $\cg$) there are no dilations $\delta_\e$. This problem can be overcome by considering a bigger metric space $\cS$ which contains $\cg$ and $\cg_\infty$ and which does possess dilating maps.

Let us consider $\cS:=\cG_{\infty}\times (0,+\infty]$ and define
\begin{eqnarray*}
\Delta_{t}: \cS &\longrightarrow& \cS \\
(v,s) &\longmapsto& (\delta_t(v),  \frac{s}{t}).\\
\end{eqnarray*}

%Let fix $W_i$  a supplementary subspace of $\cG^{(i+1)}$ in $\cG^{(i)}$ and let us consider the unique linear bijection $L: \cG \longrightarrow \cG_{\infty}$ such that:
%\begin{itemize}
%\item $L(W^i)=V^i$, where $V^i$ is the $i$-th stratum in the stratification of $\cG_{\infty}$;
%\item modulo the isomorphisms $W^i\simeq \cG^{(i)} / \cG^{(i+1)}$ and $V^i \simeq \cG^{(i)} / \cG^{(i+1)}$, $L$ corresponds to the identity map.\\
%\end{itemize}

Let us denote by $\ccG_t := \cG_{\infty}\times \{t\}$. Observe that $\ccG_1$ can be equipped with a Lie structure $[\cdot,\cdot]_1$ such that $L: \cG \longrightarrow \ccG_1$
is an isomorphism ($L$ is the identification between $\cg$ and $\cg_{\infty}$, see subsection \ref{notation}).\\
Similarly, each $\ccG_t$ can be endowed with a Lie structure  $[\cdot,\cdot]_t$ such that the dilation map
$$
\Delta_\frac{1}{t}: \ccG_1 \longrightarrow \ccG_t
$$
is an isomorphism (so we can think of $\ccG_1$ as a copy of $\cG$ in $\cS$). One can show that  $[\cdot,\cdot]_t$ depends continuously on $t\in (0,+\infty]$ (see \cite[Section 39]{Pansu}).\\
Moreover, let us denote by $G_t$ the Lie group whose corresponding Lie algebra is $\ccG_t$, by $\exp_t$ the corresponding exponential map and by
$\log_t$  its inverse.
Each $G_t$ contains a copy $\G_t$ of $\G$ which becomes denser and denser in $G_t$ as $t$ goes to infinity (see \cite[Section 39]{Pansu}).
\\

We can now define the rescaling map $h_\e$ in the following way. For each $\e>0$, let us first consider the map $\hd_{1/\e}:   G_\infty \longrightarrow G$ given by:
\begin{eqnarray*}
\xymatrix{
 G_{\infty} \ar@{->}[r]^{\log_{\infty}}  \ar@/_1.5pc/[rrrrr]_{\hd_{1/\e}}&  \cG_{\infty}  \ar@{^{(}->}[r]^{i_{1/\e}}&\ccG_{1/\e}
   \ar@{->}[r]^{{\Delta}^{-1}_{\e}} &  \ccG_1
  \ar@{->}[r]^{L}  & \cG   \ar@{->}[r]^{\exp} & G.
   }
   \end{eqnarray*}

It follows easily from the definition that
\begin{equation}\label{composizionedelta}
\hd_{\frac{1}{\a \e}}(\d_\a(\bx)) = \hd_{\frac{1}{\e}}(\bx)  \qquad \forall \; \bx\in G_\infty \; {\rm and}\; \a>0.
\end{equation}
In the following, we shall also consider the maps   $\td_{\e}:   G  \longrightarrow \cS$ (in some sense the inverses of the above maps)  given by:
\begin{eqnarray*}
\xymatrix{
G  \ar@{->}[r]^{\log}  \ar@/_1.5pc/[rrrr]_{\td_{\e}}
&
\cG \ar@{->}[r]^{L^{-1}}
 &
 \ccG_1 \ar@{->}[r]^{{\Delta}_{\e}}
 &
  \ccG_{1/\e} \ar@{^{(}->}[r]^{i} 
  & \cS.
   }
\end{eqnarray*}

Observe that $\cS$ can be endowed with a length metric -- that we shall denote by $Q$ --  which is compatible with the topology of $\cS$ and for which $\Delta_t$ are dilations (see \cite[Section 40]{Pansu}).\\

% 
% 
%  \stackrel{\log_{\infty}}{}{ \cG_{\infty}} \hookrightarrow   \cG_\e
% 
%   \xymatrix{
%    T^* G=T^*\universalcover{M} \ar@{->>}[r]\ar@{->>}[d]^{\abeliancover{\Pi}} \ar@{->>}@/_17mm/[dd]_\Pi &
%    G=\universalcover{M} \ar@{->>}[d]^{\universalcover{\pi}} \ar@{->>}@/^17mm/[dd]^\pi\\
%    T^* (\commutator{\Gamma} \backslash G)=T^* \abeliancover{M} \ar@{->>}[r] \ar@{->>}[d]^{\abeliancover{\Pi}}&
%    \commutator{\Gamma} \backslash G=\abeliancover{M} \ar@{->>}[d]^{\abeliancover{\pi}}\\
%    T^* ({\Gamma} \backslash G)=T^* {M} \ar@{->>}[r] &
%    {\Gamma} \backslash G={M}\,.
%    }
%
% 

Let us choose a compact fundamental domain $\Omega$ for the action of $\Gamma$ in $G$ (we can assume that the identity $e\in \Omega$);  for each $\e>0$ let us define the {\it rescaling maps} $h_\e : G_\infty \longrightarrow \Gamma$ such that for each $\bx\in G_\infty$, $h_\e(\bx)$ equals an element $\g\in \G$ such that $\hd_{1/\e}(\bx) \in \gamma \cdot \Omega$.\\

In particular, it follows easily from the definition of $h_\e$ and (\ref{composizionedelta}) that:
\begin{equation}\label{riscalheps}
h_{\a \e}(\d_\a(\bx)) = h_\e (\bx)  \qquad \forall \; \bx\in G_\infty \; {\rm and}\; \a>0.\\
\end{equation}

\vspace{5 pt}

\begin{Rem}
In the specific setting considered in case I, we have that $G=G_\infty$ and $\cG=\cG_\infty$ and it is easy to check that the two definitions of $h_\e$ coincide.\\
\end{Rem}

%\marginnote{\textcolor{red}{Discutere il caso degli elementi nella parte abeliana!}}

\begin{Rem}\label{propertyfe}
Properties of these maps $h_\e$:
\begin{itemize}
\item[i)]  For each $R>0$, there exists $\theta(\e)\rightarrow 0$ as $\e\rightarrow 0$ such that $h_\e$ maps the ball $B_\infty(R)$ in $G_\infty$ into $B_\G\left(\frac{R+\theta(\e)}{\e}\right)$ in $\Gamma$, which corresponds to the ball of radius $R+\theta(\e)$ with respect to the rescaled norm $\e|\cdot|$ (see \cite[\S 27]{Pansu}).
\item[ii)] If $\g_\e$ is a sequence in $\Gamma_{1/\e}$ such that $ \gamma_\e \rightarrow \bx \neq e$  as $\e\rightarrow 0$ (\ie $(\gamma_\e,\frac{1}{\e}) \rightarrow (\bx,\infty)$ in the $Q$-metric topology on $\cS$), then $\e |\gamma_\e|$ tends to $d_{\infty}(e,\bx)$ (see \cite[Proposition 41]{Pansu}).  In particular, for each $\bx\in G_\infty\setminus\{e\}$ we have that   $\td_\e(h_\e(\bx))$ converges to $\bx$; namely, if we denote by $\bx_\e:= (\log_{\infty}\bx , 1/\e)$ then
$$
Q(\td_\e (h_\e(\bx)), \bx_\e ) = \e\, Q(h_\e(\bx), \hd_{1/\e} (\bx) ) \leq \e\; {\rm diam}_Q(\Omega),
$$
which tends to $0$ as $\e$ tends to $0$.
%d_\infty(\delta_\e(\g_\e), \bx) = \e d_\infty(\g_\e, \delta_{1/\e}(\bx)) \leq C \e \,{\rm diam}(\Omega) \stackrel{\e\rightarrow 0^+}{\longrightarrow} 0.

\item[iii)] Moreover,
if $\bx \neq \by$, then $\td_\e\left( h_\e(\bx)^{-1} h_\e(\by)\right)$ tends to $\bx^{-1}\by$  and 
$\e |h_\e(\bx)^{-1} h_\e(\by)| \rightarrow d_\infty(\bx,\by)$ as $\e$ tends to zero (see  \cite[\S 27]{Pansu}); this convergence is uniform on compact sets of $G_\infty$.
\\
%%\item[iii)]
%%\item[iv)]
\end{itemize}
\end{Rem}
\end{itemize}

%%%%%%%%%%%%%%%%%%%%%%%
%%%%%%%%%%%%%%%%%%%%%%%%%%
%%%%%%%%%%%%%%%%%%%%%%%%%%%

\subsubsection{Convergence of Functions.}
Next step consists in introducing and discussing the right notion of {convergence} for functions defined on these rescaled metric spaces. We fix a reference point $x_0\in X$.\\
%. In this section we want to define and discuss what we mean by  convergence of functions with respect to these rescaled metric spaces.\\

\begin{Def}\label{defconv} Let $F_\e:  X_\e \longrightarrow \R $ and let $F : G_\infty \longrightarrow \R$. We shall say that:
\begin{itemize}
\item $F_\e$ converges pointwise to $F$, if for each $\bx\in G_\infty$ we have
$\lim_{\e \rightarrow 0^+} F_\e(h_\e(\bx)\cdot x_0) = F(\bx)$. \\
\item  $F_\e$ converges locally uniformly to $F$, if for each $R>0$ we have that
$$
\lim_{\e \rightarrow 0^+} \sup_{B_\infty(R)} |F_\e(h_\e(\bx)\cdot x_0) - F(\bx) | = 0.\\
$$
\end{itemize}
\end{Def}

\vspace{5 pt}

\begin{Rem}\label{nodependence}
\begin{itemize}
\item If $F_\e$ are equicontinuous, then the pointwise limit does not depend on the choice of $x_0$. 
In fact, if $F_\e$  are equicontinuous, then for each $\delta>0$ there exists $\eta=\eta(\delta)$ such that 
$$
d_\e(x,y)<\eta \quad \Longrightarrow \quad |F_\e(x) - F_\e(y)| <\delta.
$$
Let us now consider $x_0,x_1\in X$ and let $\bx\in G_\infty$; it is sufficient to show that for each $\delta>0$, there exists $\e_0=\e_0(\delta)$ such that
$$
|F_\e(h_\e(\bx)\cdot x_0) - F_\e(h_\e(\bx)\cdot x_1)|<\delta \qquad \forall\; 0<\e<\e_0.
$$
Observe that
\begin{eqnarray*}
d_\e( h_\e(\bx)\cdot x_0, h_\e(\bx)\cdot x_1) &=&
\e d ( h_\e(\bx)\cdot x_0, h_\e(\bx)\cdot x_1) 
= \e d(x_0, x_1),
\end{eqnarray*}
where we used the $\G$-left-invariance of $d$. Hence, it follows from the equicontinuity of $F_\e$, that it is sufficient to 
choose $\e_0 < d(x_0, x_1)/\eta$.\\

\item Similarly, if $F_\e$ are equicontinuous, then the locally uniform limit does not depend on the choice of $x_0$. In fact, if $x_0,x_1\in  X$ it follows from what we have  said above that for each $\delta>0$
$$
\lim_{\e \rightarrow 0^+} \sup_{B_\infty(R)} |F_\e(h_\e(\bx)\cdot x_0) - F(h_\e(\bx)\cdot x_1) | \leq \delta
$$
and consequently (for the arbitrariness of $\delta$)
$$
\lim_{\e \rightarrow 0^+} \sup_{B_\infty(R)} |F_\e(h_\e(\bx)\cdot x_0) - F(h_\e(\bx)\cdot x_1) | =0.
$$
Hence, using the triangle inequality (and the fact that the sup of a sum is less or equal than the sum of the sup's):
$$
\lim_{\e \rightarrow 0^+} \sup_{B_\infty(R)} |F_\e(h_\e(\bx)\cdot x_0) - F(\bx)|=0
\quad 
\Longrightarrow 
\quad
\lim_{\e \rightarrow 0^+} \sup_{B_\infty(R)} |F_\e(h_\e(\bx)\cdot x_1) - F(\bx)) | =0.\\
$$
%\item  If $F_\e$ are equiLipschitz on compact sets, then the above limits do not depend on the choice of $x_0$. In fact, if $x_1\in X$ then for $\bx \in B_\infty(R)$:
%\begin{eqnarray*}
%\left| F_\e( h_\e(\bx)\cdot x_0) - F_\e(  h_\e(\bx)\cdot x_1) |\right| &\leq& K d_\e \left( h_\e(\bx)\cdot x_0), h_\e(\bx)\cdot x_1)  \right) = \\
%&=& \e K d \left( h_\e(\bx)\cdot x_0), h_\e(\bx)\cdot x_1)  \right) = \\
%&=& \e K d ( x_0, x_1)  \stackrel{\e\rightarrow 0^+}{\longrightarrow 0},
%\end{eqnarray*}
%\noindent where we have used that $d$ is $\G$-invariant and that $h_\e(B_\infty(R))$ are all contained in a bounded ball  of $\Gamma$ and therefore $\{h_\e(B_\infty(R))\cdot x_0\} \cup \{h_\e(B_\infty(R))\cdot x_1\}$ Lie in compact set of $X$ and therefore one can find a common Lipschitz constant $K=K(R, x_0, x_1)$.\\
\end{itemize}
\end{Rem}

\vspace{10 pt}

\vspace{10 pt}

\subsection{{Rescaled Ma\~n\'e Potentials and Generalized Mather's $\beta$-function}}\label{secgenbeta}

Let $T>0$; for each $0<\e<1$ consider  the {\it rescaled Ma\~n\'e potential} $\phi^\e$, {\it i.e.}, the Ma\~n\'e potential associated to the rescaled metric $d_\e$. Observe that
\begin{eqnarray*}
\phi^\e(y,x,S) &:= &\inf \left\{ \int_0^S L (\eta(s),\e \dot{\eta}(s))\,ds\; \Big| \; \eta\in C^1([0,S], X),\; \eta(0)=y, \;\eta(S)=x \right\} \\
&=&
\e \inf \left\{ \int_0^{\frac{S}{\e}} L (\eta(t),\dot{\eta}(t))\,ds\; \Big| \; \eta\in C^1([0,S/\e], X),\; \eta(0)=y, \;\eta(S/\e)=x \right\} \\
&=& \e \phi(y,x,S/\e),\\
\end{eqnarray*}
where $\phi$ denotes the standard Ma\~n\'e potential associated to $L$ with the metric $d$, as defined in  (\ref{manepotential}) (see also \cite{Gonzalobook, SorrentinoLectureNotes} for more details).\\

Recalling the definition of $h_\e$ (see subsection \ref{subsecheps}), let us define the following family of functions ($x_0$ is a fixed reference point in $X$): 
\begin{eqnarray}\label{Feps}
F_\e: G_\infty \times (0,+\infty) &\longrightarrow& \R\nonumber\\
(\bx,T) &\longmapsto& \phi^\e(x_0, h_\e(\bx)\cdot x_0, T). \\
\nonumber
\end{eqnarray}

We want to discuss some properties of these functions.\\

\begin{Prop}\label{equiLipschitzphi}
Let $\lambda>0$ and $\K$ be a compact subset of $G_\infty$. Then, for $\e$ sufficiently small, the maps $F_\e$ are equiLipschitz on $\K \times [\lambda, +\infty)$.\\
\end{Prop}

\begin{Proof}
Let $\l>0$ and
let us prove the statement for closed balls of radius $\lambda/8$ in $(G_\infty,d_\infty)$ (the proof for a general compact set follows by covering it with a finite number of such balls\footnote{Observe that as $\lambda$ approaches zero, then the Lipschitz constant tends to infinity since a larger and larger number of balls is needed to cover the space.}).
Observe that if 
$\bx,\by$ belong to such a ball, then it follows from Remark \ref{propertyfe} (iii) that for $\e$ sufficiently small
$$
\e\, d(h_\e(\bx)\cdot x_0, h_\e(\by)\cdot x_0) < 2\, d_\infty(\bx, \by) \leq 2 \cdot \frac{\l}{4} =\frac{\l}{2}.
$$ 
Therefore, if $T\geq \lambda$, then $\frac{T}{2\e} > \frac{\l}{2}$ and  using Lemma \ref{lemma1} with $K=K(\l)$ and Remark \ref{propertyfe} (iii), we obtain:
\begin{eqnarray*}
|F_\e(\bx,T) - F_\e(\by,T)| &=& 
 |\phi^\e(x_0, h_\e(\bx)\cdot x_0, T) - \phi^\e(x_0, h_\e(\by)\cdot x_0, T)| \\
&=& \e\, |\phi (x_0, h_\e(\bx)\cdot x_0, T/\e) - \phi (x_0, h_\e(\by)\cdot x_0, T/\e)| \\
&\leq& \e K d ( h_\e(\bx)\cdot x_0,  h_\e(\by)\cdot x_0)  \\
&=& \e K d_{\G,x_0} ( h_\e(\bx),  h_\e(\by))  
\leq \e c_1 K | h_\e(\bx)^{-1}  h_\e(\by)| \\
&\leq&  2\,c_1\, K d_\infty(\bx,\by)
 =: K_1d_\infty(\bx,\by),
\end{eqnarray*}
where $c_1$ is a constant (independent of $\lambda$) which relates the orbit metric to the word metric (a-priori the norm that we are considering on $\G$ might come from the word metric for some set of generators); see subsection \ref{metrics} and Proposition \ref{prop1}.

To prove that $F_\e$ are Lipschitz in the $T$-variable, let us observe that the Ma\~n\'e potential is Lipschitz in the time component for $T\geq \lambda>0$ (the Lipschitz constant depends on $\l$; see \cite[Proposition 3-4.1]{Gonzalobook}). Then, if we denote by $C=C(\l)$ this Lipschitz constant, we obtain for $S,T\geq \l$:
\begin{eqnarray*}
|F_\e(\bx,T) - F_\e(\bx,S)| &=& 
 |\phi^\e(x_0, h_\e(\bx)\cdot x_0, T) - \phi^\e(x_0, h_\e(\bx)\cdot x_0, S)| \\
&=& \e\, |\phi (x_0, h_\e(\bx)\cdot x_0, T/\e) - \phi (x_0, h_\e(\bx)\cdot x_0, S/\e)| \\
&\leq& \e C \left|\frac{T}{\e}-\frac{S}{\e} \right|  
= C \,|T-S|.
\end{eqnarray*}

% we have used Lemma \ref{lemma1} and . In fact, if $\e$ is sufficiently small, then
%$$
%d(h_\e(\bx), h_\e(\by)) = \frac{1}{\e} d_\e(h_\e(\bx), h_\e(\by)) \leq \frac{C}{\e} d_\infty(\bx,\by);
%$$
%hence it follows from Lemma \ref{lemma1} that for $\e$ small that $\Phi(x_0, h_\e(\bx)\cdot x_0, T/\e$ is Lipschitz on the set $\{T > d_\infty(\bx,\by)\}$ and the Lipschitz constant is independent of $\e$.\\
\end{Proof}

%%%%%%%%%%%%
Let us now prove the following rescaling property.

\begin{Prop} \label{rescalingprop}
 For each $\e>0$, $\bx\in G_\infty$ and $T,S>0$, the following rescaling property holds:
$$
F_\e(\bx, S) = \frac{S}{T}\, F_{\e \frac{T}{S}} \Big(   
\delta_{\frac{T}{S}}(\bx), T
\Big).\\
$$
\end{Prop}

\begin{Proof}
It suffices to consider the definition of $F_\e$ and $\phi^\e$, and to apply the formula of change of variables in the integral. In fact, if $\gamma: [0,S]\longrightarrow X$ is a curve with end-points $x_0$ and $h_\e(\bx)\cdot x_0$, then one can reparametrize it to obtain a curve $\eta: [0,T]\longrightarrow X$ with the same end-points (\ie $\eta(t)=\gamma(\frac{S}{T}t)$) and such that
$$
\int_0^S L(\g(s), \e \dot{\g}(s)) \,ds =  \frac{S}{T} \int_0^T L\left(\eta(t),  \frac{\e T}{S} \dot{\eta}(t)\right) \,dt.
$$
Therefore, 
\begin{eqnarray*}
F_\e(\bx, S) &=&
\phi^\e(x_0, h_\e(\bx)\cdot x_0,S) 
= \frac{S}{T} \phi^{\e\frac{T}{S}}(x_0, h_\e(\bx)\cdot x_0,T) \nonumber\\
&=&  \frac{S}{T} \phi^{\e\frac{T}{S}} \Big(
x_0, h_{\e \frac{T}{S}} \left( \delta_{\frac{T}{S}}(\bx)  \right)\cdot x_0,T
\Big) 
=  \frac{S}{T}\, F_{\e \frac{T}{S}} \Big(   
\delta_{\frac{T}{S}}(\bx), T
\Big),
\end{eqnarray*}
where in the second-last equality we used (\ref{riscalheps}).\\
\end{Proof}

\subsubsection{Generalized Mather's $\beta$-Function}

For each $T>0$ and $\bx \in G_\infty$ we can define the following function: 
\begin{eqnarray}\label{convgenerbetafunction}
\bb(\bx, T) :=  \lim_{\e\rightarrow 0^+} F_\e(\bx, T) =
 \lim_{\e\rightarrow 0^+} \phi^\e(x_0, h_\e(\bx)\cdot x_0, T).
\end{eqnarray}
The existence of this limit will be discussed in Proposition \ref{whatislimit}. For now, let us assume that this limit exists and let us prove the main properties of this {\it generalized} ({\it time-$T$}) {\it  Mather's $\beta$ function}  $\bb(\cdot, T)$.\\

First of all, observe that this limit  does not depend on the particular choice of $x_0$ since $\phi^\e$ are locally equiLipschitz (with respect to $d_\e$), hence equicontinuous (see Remark \ref{nodependence}). Moreover, it enjoys a rescaling property similar to what we have seen in Proposition \ref{rescalingprop}.\\

%\begin{Def}[{\bf Generalised Mather's $\beta$-function}] 
%We define the Generalised Mather's $\beta$-function as 
%\begin{eqnarray*}
%\bbb : G &\longrightarrow& \R\\
%\bx  &\longmapsto& \bb(\bx,1).\\
%%
%%\bbb : \cG &\longrightarrow& \R\\
%%v &\longmapsto& \bb(\exp(v),1),
%%
%\end{eqnarray*}
%%where $\exp: \cG \longrightarrow G$ denotes the (Lie) exponential map.\\
%\end{Def}

%\begin{Rem}
%\begin{itemize}
%%\item One could also consider the function $\bbb$ directly on $G$ rather than on its Lie Algebra.
%\item Observe that if $X$ is the abelian cover of a compact smooth manifold $M$ without boundary (see examples in subsection \ref{sec1.1}), then $\Gamma \simeq \Z^k$ (
%$k=b_1(M)$ is the first Betti number of $M$) and $G=H_1(M;\R)\simeq \R^k$; in this case $\bbb$ is exactly Mather's $\beta$-function (see \cite[Proposition 1]{Mather91}).
%\item If $X$ is the universal cover of a compact smooth manifold $M$ without boundary, then $\Gamma=\pi_1(M)$ and $\bbb$ is a sort of {\it homotopical} version of Mather's $\beta$-function, which extends this function to the Lie algebra of the asymptotic cone of the fundamental group.\\
%\end{itemize}
%\end{Rem}

%Let us discuss some properties of $\bb(\cdot, T)$. \\%and $\bbb(\cdot)$.\\

\begin{Prop}\label{riscalamentobeta}
For each $T,S>0$ and for each $\bx\in G_\infty$, we have:
$$
\frac{1}{S} \bb(\bx, S) = \frac{1}{T} \bb(\delta_{T/S}(\bx), T).
$$
In particular, for each $T>0$ and $\bx \in G_\infty$ we have:
$$
\bb(\bx, T) = {T} \bb\left(\delta_{1/T}(\bx), 1\right).\\% =  {T} \bbb\left( \delta_{1/T}(\bx)\right).\\
$$
%where   $\delta_s$ denote the dilations seen as automorphisms of  $G$.\\
\end{Prop}

\vspace{5 pt}

\begin{Proof} It is an easy consequence of Proposition \ref{rescalingprop}:
\begin{eqnarray*}
 \bb(\bx, S) &=& \lim_{\e\rightarrow 0^+}  F_\e(\bx, S) \\
 &=& \lim_{\e\rightarrow 0^+} \left(\frac{S}{T}\, F_{\e \frac{T}{S}} \Big(   
\delta_{\frac{T}{S}}(\bx), T
\Big) \right) 
= \frac{S}{T}  \lim_{\e\rightarrow 0^+}  F_{\e \frac{T}{S}} \Big(   
\delta_{\frac{T}{S}}(\bx), T
\Big)  \\
&=&\frac{S}{T}  \lim_{\e'\rightarrow 0}  F_{\e'} \Big(   
\delta_{\frac{T}{S}}(\bx), T
\Big)  
= \frac{S}{T} \bb(\delta_{T/S}(\bx), T).
\end{eqnarray*}

%We know that $\bb(\bx, S) =  \lim_{\e \rightarrow 0^+} \phi^\e(x_0, h_\e(\bx)\cdot x_0,S)$.
%Observe that:
%\begin{eqnarray}\label{cambiovar}
%\phi^\e(x_0, h_\e(\bx)\cdot x_0,S) &=& \frac{S}{T} \phi^{\e\frac{T}{S}}(x_0, h_\e(\bx)\cdot x_0,T) =\nonumber\\
%&=&  \frac{S}{T} \phi^{\e\frac{T}{S}} \Big(
%x_0, h_{\e \frac{T}{S}} \left( \delta_{\frac{T}{S}}(\bx)  \right)\cdot x_0,T
%\Big).
%\end{eqnarray}
%It is sufficient, in fact, to recall the definition of $\phi^\e$ and to observe that any curve $\sigma(s)$ on $[0,S]$  can be reparameterized to a curve $\eta(t)$ on $[0,T]$, simply by taking $s=t \frac{S}{T}$ (and vice-versa); in particular,  $ds= \frac{S}{T} dt$.\\ Taking now the limit as $\e$ goes to zero on both sides we obtain:
%$$
%\bb(\bx, S) \longleftarrow \phi^\e(x_0, h_\e(\bx)\cdot x_0,S) = \frac{S}{T} \phi^{\e\frac{T}{S}} \Big(
%x_0, h_{\e \frac{T}{S}} \left( \delta_{\frac{T}{S}}(\bx)  \right)\cdot x_0,T\Big)
%\longrightarrow
% \frac{S}{T} \bb\left( \delta_{\frac{T}{S}}(\bx), T\right),
%$$
%which concludes the proof.\\
\end{Proof}

Let us prove other two important features of $\bb(\cdot, 1)$.

\begin{Prop}\label{propfunzionebeta}
{\rm I.} $\bb(\cdot, 1) : G_\infty \longrightarrow \R$ is a convex function, {\it i.e.}, for each $\lambda\in (0,1)$ and for each $\bx, \by \in G_\infty$ we have:
$$
\bb\Big(\delta_\l(\bx) \cdot \delta_{1-\l}(\by),1\Big) \leq \lambda \bb(\bx,1) + (1-\l) \bb(\by,1).
$$
{\rm II.} 
$\bb(\cdot,1)$ is a superlinear function, {\it i.e.}, for each $A>0$ there exists $B=B(A)\geq 0$ such that
$$
\bb (\bx,1) \geq A\; d_{\infty}(e,\bx) - B.\\
$$
\end{Prop}

%\textcolor{red}{(se fosse vera la definizione con il limite, probabilmente si potrebbe definire anche la convessit\`a.)}

\vspace{10 pt}

\begin{Proof}
%I. Let us the following notation: $\bx_v := \exp^{-1}(v)$;
Let us start by observing that:
\begin{eqnarray}\label{triangolare}
\phi^\e (x_0, h_\e(\delta_\l(\bx) \cdot \delta_{1-\l}(\by)) \cdot x_0, 1)
& \leq &
\phi^\e (x_0, h_\e(\delta_\l(\bx))\cdot x_0, \lb)  \nonumber\\
&&+\; \phi^\e (h_\e(\delta_\l(\bx))\cdot x_0, h_\e(\delta_\l(\bx) \cdot \delta_{1-\l}(\by))\cdot x_0, 1-\lb).
\end{eqnarray}
%
%
%
%Let us start by observing that:
%\begin{eqnarray}\label{triangolare}
%\phi^\e (x_0, h_\e(\bx_{\delta_\l(v) + \delta_{1-\l}(w)})\cdot x_0, 1)
%& \leq &
%\phi^\e (x_0, h_\e(\bx_{\delta_\l(v)})\cdot x_0, \lb) + \nonumber\\
%&&+\; \phi^\e (h_\e(\bx_{\delta_\l(v)})\cdot x_0, h_\e(\bx_{\delta_\l(v)+ \delta_{1-\l}(w)})\cdot x_0, 1-\lb).
% %
%\end{eqnarray}

%
Using Proposition  \ref{rescalingprop} with $T=1$ and $S=\lambda$, % and recalling the relation between dilations in $G$, dilations in $\cG$ and the exponential map), 
we obtain that:
\begin{eqnarray}\label{riscalamento1}
\phi^\e (x_0, h_\e(\delta_\l(\bx))\cdot x_0, \lb) &=& 
\lambda \phi^{\e/\l} (x_0, h_{\e/\l}(\delta_{1/\lambda}\delta_\l(\bx))\cdot x_0, 1)  \nonumber\\
%&=& \lambda \phi^{\e/\l} (x_0, h_{\e/\l}(\bx_{\delta_{1/\l}(\delta_\l(v))})\cdot x_0, 1) =\nonumber\\
&=& \lambda \phi^{\e/\l} (x_0, h_{\e/\l}(\bx)\cdot x_0, 1).
\end{eqnarray}

%Let us now denote $\g_1:= h_\e(\bx_{\delta_\l(v)})$ and $\g_2:=h_\e(\bx_{\delta_\l(v)+ \delta_{1-\l}(w)})$. 
It follows from Remark \ref{propertyfe} (iii) that for every $\bx\neq \by$, one has
$$
d_\e\Big(\left( h_\e (\bx) \right)^{-1} h_\e (\by), \; h_\e(\bx^{-1}\by )  \Big) \longrightarrow 0 \quad {\rm as}\; \e\rightarrow 0.
$$
Using the local equiLipschitzianeity of $F_\e$ (see Proposition \ref{rescalingprop}) and the fact that $L$ is invariant under the action of $\Gamma$, we know that (similarly to what done in (\ref{riscalamento1})):
\begin{eqnarray}\label{riscalamento2}
&& \phi^\e (h_\e(\delta_\l(\bx))\cdot x_0, h_\e(\delta_\l(\bx) \cdot \delta_{1-\l}(\by))\cdot x_0, 1-\lb) \nonumber \\
&&\qquad = \; 
\phi^\e \left( x_0,  \left(h_\e \big(\delta_\l (\bx) \right)^{-1} h_\e \big(\delta_\l(\bx) \cdot \delta_{1-\l}(\by)\big) \cdot x_0, 1-\l  \right) \nonumber\\
&&\qquad = \; 
\phi^\e \left( x_0,  h_\e( \delta_\l (\bx)^{-1}\delta_\l(\bx)  \delta_{1-\l}(\by)) \cdot x_0, 1-\l  \right) + o_\e(1)\\
&&\qquad = \; 
\phi^\e \left( x_0,  h_\e(  \delta_{1-\l}(\by)) \cdot x_0, 1-\l  \right) + o_\e(1) \nonumber\\
%&=& \phi^\e \left( x_0,  h_\e(\bx_{\delta_{1-\l}(w)}) \cdot x_0, 1-\l  \right) =\nonumber\\
&&\qquad = \;  (1-\lambda) \phi^{\frac{\e}{1-\l}} ( x_0,  h_{\frac{\e}{1-\l}}(\by) \cdot x_0, 1  ) + o_\e(1). \nonumber
\end{eqnarray}

Taking the limit as $\e$ goes to zero in (\ref{triangolare}), using (\ref{riscalamento1}), (\ref{riscalamento2}) and the definition of $\bb$, we conclude that
$$
\bb\Big(  \d_\l(\bx) \cdot \d_{1-\l}(\by) ,1 \Big) \leq \lambda \bb(\bx,1) + (1-\l) \bb(\by,1).
$$
%and therefore (see definition of $\bbb$):
%$$
%\bbb\Big(  \d_\l(\bx) \cdot \d_{1-\l}(\by)\Big) \leq \lambda \bbb(\bx) + (1-\l) \bbb(\by).\\
%$$

\noindent II. 
 Since $L$ is superlinear, then for each $A>0$ there exists $B=B(A)\geq0$ such that for each $(x,v)\in TX$
$$
L(x,v) \geq A \|v\|_x - B,
$$
where $\|\cdot\|_x$ is the norm associated to the $\Gamma$-invariant Riemannian metric that we are considering on $X$. Then, for each 
$\eta: [0,T] \rightarrow X$ joining $x$ to $y$ we have:
$$
\int_0^T L(\eta(t),\dot{\eta}(t))\,dt \geq \int_0^T A \|\dot{\eta}(t)\| \,dt - BT \geq A\, d(x,y) - BT.
$$
Then, it follows easily that
\begin{eqnarray*}
\phi^\e(x_0, h_\e(\bx)\cdot x_0, 1) &\geq& \e A\, d\left(x_0, h_\e(\bx)\cdot x_0
\right) - B 
=
A\, \e\,d_{\G,x_0}\left(e, h_\e(\bx)
\right) - B.
\end{eqnarray*}
Taking the limit as $\e$ goes to zero (recall Remark \ref{propertyfe}) we can conclude that:
$$
\bb(\bx,1) \geq  A \,d_{\infty}(e, \bx) - B.
$$
\end{Proof}

\vspace{10 pt}

\begin{Rem}\label{limsupbeta}
Observe that for proving Propositions \ref{riscalamentobeta} and \ref{propfunzionebeta} it would be sufficient to consider the limsup in (\ref{convgenerbetafunction}), rather than the limit (however, this is not enough for the proof of the main theorem). Therefore,
if we denoted by 
$$\bb^+(\bx, T) := \limsup_{\e\rightarrow 0^+} F_\e(\bx, T) =
 \limsup_{\e\rightarrow 0^+} \phi^\e(x_0, h_\e(\bx)\cdot x_0, T),$$
where $T>0$ and $\bx\in G_\infty$, then this function would satisfy the same rescaling property and would remain convex (in the same sense as in Proposition \ref{propfunzionebeta}) and superlinear.\\
On the other hand, if one considers
$$\bb^-(\bx, 1) = \liminf_{\e\rightarrow 0^+} \phi^\e(x_0, h_\e(\bx)\cdot x_0, 1)$$
then this function  would continue to satisfy the rescaling property in Proposition \ref{riscalamentobeta}, but a-priori convexity and superlinearity could not be deduced from the same arguments as in Proposition \ref{propfunzionebeta}.\\
\end{Rem}

We can now define one of the main object in the statement (and in the proof) of the main result (see statement in subsection \ref{secmainres}).

\begin{Def} We define the {\it Generalized Mather's $\beta$-function} as
\begin{eqnarray*}
\bL:   G_\infty  &\longrightarrow& \R\\
\bx &\longmapsto& \bb(\bx,1).
\end{eqnarray*}
\end{Def}

\vspace{10 pt}

We can now prove the following result, which describes what is $\bL$ and therefore -- using Proposition \ref{riscalamentobeta} -- proves the existence of the limit in (\ref{convgenerbetafunction}). Let us first recall that an absolutely continuous curve
$\g: [a,b] \longrightarrow G_\infty$ is said to be {\it horizontal} if  for almost every $t\in [a,b]$ the tangent vector $\dot{\gamma}(t)\in V_1$, \ie it lies in the first layer of the graded algebra $\fg_{\infty}$.\\
\begin{Prop}\label{whatislimit}
Let $\beta: H_1( X\slash\Gamma; \R)\longrightarrow \R$ be Mather's $\beta$-function associated to the projection of $L$ to $T( X\slash\Gamma)$ and let us consider the projection 
$\opi: \cg_\infty \longrightarrow \frac{\cG_\infty}{[\cG_\infty,\cG_\infty]} \simeq \frac{\cG}{[\cG,\cG]} \simeq  H_1( X\slash\Gamma; \R)$. Then, for each $\bx \in G_\infty$ we have:
$$ \bL(\bx):=\lim_{\e\rightarrow 0^+} \phi^\e(x_0, h_\e(\bx)\cdot x_0, 1) =  \inf_{\sigma \in \cH_{\bx}} \int_0^1 \beta(\opi(\dot{\sigma}(s)))\,ds,
$$
where $\cH_{\bx}$ denotes the set of  horizontal curves $\sigma: [0,1] \longrightarrow G_\infty$ connecting $e$ to $\bx$.\\
\end{Prop}

\begin{Rem}
({\it i}). Note that if $\Gamma=\Z^k$, then $G_\infty=\cg_{\infty}=\R^k$ and it follows from the convexity of $\beta$ and Jensen's inequality that $\inf_{\sigma \in \cH_{\bx}} \int_0^1 \beta(\opi(\dot{\sigma}(s)))\,ds = \beta(\bx)$. This is consistent with Mather's result  \cite[Proposition 1]{Mather91}.\\
({\it ii}) Observe that  this expression of $\bL(\bx)=\bb(\bx,1)$ coincides with the {\it value function} $\cL(\cdot,1)$ introduced in \cite[section 2.3]{BCP}. Hence, all properties of $\bb$ discussed in Propositions 
\ref{riscalamentobeta} and \ref{propfunzionebeta} could be also deduced from analogous properties of $\cL$ proved in \cite[Proposition 2.4]{BCP}.\\
\end{Rem}

Before proving Proposition \ref{whatislimit}, let us discuss what happens when $\bx \in \exp_{\infty}(V_1)$, \ie what we could call the ``abelian part''. We shall prove that in this case the limit in (\ref{convgenerbetafunction}) exists and is related to Mather's $\beta$-function. Hereafter, we shall identify  $H_1(X\slash\Gamma; \R) \simeq \frac{\cG_\infty}{[\cG_\infty,\cG_\infty]} $.\\

\begin{Rem}
By assumption we are considering the case in which  $H_1(X;\R)=0$; otherwise,  the proof is essentially the same, but  since the map $\opi$ is not surjective,  $\frac{\cg}{[\cg,\cg]}$ will be identified with a (possibly) proper subgroup of $H_1( X\slash\Gamma;\R)$. This corresponds, for example, to the abelian subcover case discussed in \cite{CIS}. \\
\end{Rem}
 
Let us consider $\widehat{L}: T(\widehat{ X\slash\Gamma}) \longrightarrow \R$ the associated Tonelli Lagrangian on the maximal free abelian cover of $ X\slash\Gamma$ (see, for example, subsection \ref{subsecexamples}). Let $\bp: X \rightarrow \widehat{ X\slash\Gamma}$ be the projection on this cover and let us denote by $\widehat{\phi}$ the associated Ma\~n\'e potential. It is easy to check that:
\begin{equation}\label{sollevamentopotenziale}
\phi(x,y,T) \geq \widehat{\phi}(\bp(x), \bp(y), T)
\end{equation}
(essentially, it follows from the definition, since $\widehat{\phi}$ is obtained by taking the infimum over a set of curves, which is larger than the projected set of curves in $X$ connecting $x$ to $y$).

Hence:
\begin{eqnarray}\label{eq13}
\phi^\e(x_0, h_\e(\bx)\cdot x_0, 1) &=& \e \phi (x_0, h_\e(\bx)\cdot x_0, 1/\e) \nonumber\\
&\geq& \e \widehat{\phi} (\bp(x_0), \bp(h_\e(\bx)\cdot x_0), 1/\e) \stackrel{\e\rightarrow 0^+}{\longrightarrow} \beta(\opi(\log_{\infty} \bx)),
\end{eqnarray}
where in the last step we have used  \cite[Proposition 1]{Mather91}. \\

\begin{Lem}\label{lemmacasoabeliano} 
If $\bx \in \exp_\infty(V_1)\leq G_\infty$, then 
$$ \lim_{\e \rightarrow 0^+}\phi^\e(x_0, h_\e(\bx)\cdot x_0, 1) = \beta(\opi(\log_\infty \bx)). \\$$
In particular, $\bL(\bx):=\bb(\bx,1)=\beta(\opi(\log_\infty \bx)).$
\end{Lem}
\begin{Proof}
In the light of (\ref{eq13}), it is sufficient to prove that $\limsup_{\e \rightarrow 0^+}\phi^\e(x_0, h_\e(\bx)\cdot x_0, 1) \leq \beta(\opi(\log_\infty \bx))$.
Let us consider the projection $\pi_\Gamma: \Gamma \longrightarrow \G\cap \exp(V_1)$ (recall that we are identifying $W_1\simeq V_1$), defined in the following way.
Following Malcev \cite{Malcev} (see also \cite[section 36]{Pansu}), the subgroup of commutators $[\G,\G]$ coincides with $\G\cap [G,G]$, and the image $\pi(\Gamma)$ is a co-compact lattice in $V_1$ ($\pi$ denotes here the projection from $G$ to $\frac{\cg}{[\cg,\cg]}$). Then:

\begin{eqnarray*}
\xymatrix{
\G  \ar@{->}[r]  \ar@/_1.5pc/[rr]_{\pi_{\G}}
&
\pi(\Gamma) \leq \frac{\cg}{[\cg,\cg]} \simeq
V_1 \ar@{->}[r]^{{\exp}}
 &
\Gamma \cap \exp(V_1)
}
\end{eqnarray*}

Then, using the Lipschitzianeity of $\phi$:
\begin{eqnarray*}
\phi^\e(x_0, h_\e(\bx)\cdot x_0, 1) &=& \e \phi (x_0, h_\e(\bx)\cdot x_0, 1/\e) \\
&\leq& \e \phi (x_0, \pi_\G(h_\e(\bx))\cdot x_0, 1/\e) + C  d_\e(\pi_\G(h_\e(\bx)) \cdot x_0, h_\e(\bx)\cdot x_0)\\
&\leq& \e \widehat{\phi} \left(\bp(x_0), \bp\left(\pi_\G(h_\e(\bx))\cdot x_0\right), 1/\e\right) + C' \e |\pi_\G(h_\e(\bx))^{-1} h_\e(\bx)| \\
& \stackrel{\e\rightarrow 0^+}{\longrightarrow}& \beta(\opi(\log \bx)),
\end{eqnarray*}
where in the last step we have used that $ \e |\pi_\G(h_\e(\bx))^{-1} h_\e(\bx)|$ tends to zero (since $\bx\in \exp(V_1)$) and the fact that 
$\phi (x_0, \pi_\G(h_\e(\bx))\cdot x_0, 1/\e) = \widehat{\phi} (\bp(x_0), \bp\left(\pi_\G(h_\e(\bx))\cdot x_0\right), 1/\e)$, where the equality in (\ref{sollevamentopotenziale}) holds because  the minimizer in the abelian cover lifts to a curve connecting  $x_0$ to $\pi_\G(h_\e(\bx))\cdot x_0$   and consequently this lifted curve must be the minimizer in $X$ as well.
\end{Proof}

We can now prove Proposition \ref{whatislimit}.

\begin{Proof}[{\bf Proposition \ref{whatislimit}}]
As observed in Remark \ref{limsupbeta}, if $\bb^+(\bx, T)$ denotes the $\limsup_{\e\rightarrow 0^+} \phi^\e(x_0, h_\e(\bx)\cdot x_0, T)$, then this function satisfies the same rescaling and convexity properties as $\bb$.\\
Let $w:[0,1] \longrightarrow G_\infty$ be a horizontal curve connecting $e$ to $\bx$ and let $N\in \N$. We define $\bx^N_k:=w(\frac{k}{N}) \in G_\infty$, for $k=0,\ldots, N$ (where $\bx^N_0=e$ and $\bx^N_N=\bx$). Then:
\begin{eqnarray*}
\phi^\e(x_0, h_\e(\bx)\cdot x_0, 1)& \leq& \sum_{k=0}^{N-1} \phi^\e\left(h_\e(\bx^N_k)\cdot x_0, h_\e(\bx^N_{k+1})\cdot x_0, \frac{1}{N}\right) 
=  \sum_{k=0}^{N-1} \phi^\e\left(x_0, (h_\e(\bx^N_k)^{-1}h_\e(\bx^N_{k+1}))\cdot x_0, \frac{1}{N}\right).
\end{eqnarray*}

Therefore, using the fact that $w$ is horizontal and Lemma \ref{lemmacasoabeliano}:
\begin{eqnarray*}
\bb^+(\bx,1)&=& \limsup_{\e \rightarrow 0^+} \phi^\e(x_0, h_\e(\bx)\cdot x_0, 1) \\
&\leq&  \sum_{k=0}^{N-1} \limsup_{\e\rightarrow 0^+}  \phi^\e\big(x_0, h_\e((\bx^N_k)^{-1} \bx^N_{k+1})\cdot x_0, \frac{1}{N}\big) \\
&=&  \sum_{k=0}^{N-1} \bb^+ \left((\bx^N_k)^{-1}\bx^N_{k+1}, \frac{1}{N}\right) 
=  \sum_{k=0}^{N-1} \frac{1}{N} \bb^+ \left(\delta_N((\bx^N_k)^{-1}\bx^N_{k+1}), 1\right) \\
&=&  \sum_{k=0}^{N-1} \frac{1}{N} \b \left( \opi \log_{\infty} (\delta_N((\bx^N_k)^{-1}\bx^N_{k+1})) \right) \stackrel{N\rightarrow +\infty}{\longrightarrow} \int_0^1 \beta(\opi(\dot{w}(s)))\,ds.
\end{eqnarray*}

Hence, 
$$\bb^+(\bx,1)
\leq  \inf_{w \in \cH_{\bx}} \int_0^1 \beta(\opi(\dot{w}(s)))\,ds.\\
$$

\vspace{20 pt}

Let us now denote by 
$\bb^-(\bx, 1)$  the liminf in (\ref{convgenerbetafunction}), \ie
$\bb^-(\bx, 1) = \liminf_{\e\rightarrow 0^+} \phi^\e(x_0, h_\e(\bx)\cdot x_0, 1)$ (see Remark \ref{limsupbeta}) and let us prove that
$$
\bb^-(\bx,1) \geq  \inf_{\sigma \in \cH_{\bx}} \int_0^1 \beta(\opi(\dot{\sigma}(s)))\,ds.\\
$$

Let us consider a subsequence $\e_j\rightarrow 0^+$ such that
$$
\lim_{\e_j\rightarrow 0^+} \phi^{\e_j}(x_0, h_{\e_j}(\bx)\cdot x_0, 1) =
\liminf_{\e\rightarrow 0^+} \phi^\e(x_0, h_\e(\bx)\cdot x_0, 1) = \bb^-(\bx,1).
$$
For every $j$, let $\s_j:[0, \frac{1}{\e_j}] \longrightarrow X$ be the minimizer connecting $x_0$ to $h_{\e_j}(\bx)\cdot x_0$; in particular, it follows from the minimization property that for every $N\in \N$:
$$
 \phi^{\e_j}(x_0, h_{\e_j}(\bx)\cdot x_0, 1) = \sum_{k=0}^{N-1}  \phi^{\e_j}\left(   \s_j\left(\frac{k}{N}\right), \s_j\left(\frac{k+1}{N}\right), \frac{1}{N}   \right). 
$$
Since $\Gamma \cdot x_0$ becomes denser and denser in $X_{\e_j}$ (it follows from the fact that the action is co-compact), then we can find elements $\g^{(j,N)}_k \in \Gamma$ such that
$ d_{\e_j}(\g^{(j,N)}_k \cdot x_0, \s_j (k/N)) \longrightarrow 0$ as $\e_j\rightarrow 0$ (for each $k=0, \ldots, N$). Observe that we can choose $\g^{(j,N)}_0=e$ and $\g^{(j,N)}_N=h_{\e_j}(\bx)$.
In particular, up to extracting a common subsequence, we can assume that for $k=0, \ldots, N$,
$\g^{(j,N)}_k \longrightarrow \bar{x}^N_k \in G_\infty$ as $\e_j\rightarrow 0$.
Therefore, using the Lipschitzianeity of $\phi^{\e_j}$ (see Proposition \ref{equiLipschitzphi}):
\begin{eqnarray*}
\phi^{\e_j}(x_0, h_{\e_j}(\bx)\cdot x_0, 1) &= & \sum_{k=0}^{N-1}  \phi^{\e_j} \left(\s_j\left(\frac{k}{N}\right), \s_j\left(\frac{k+1}{N}\right), \frac{1}{N} \right) \\
&=& \sum_{k=0}^{N-1}  \phi^{\e_j}\left(\g^{(j,N)}_k \cdot x_0, \g^{(j,N)}_{k+1} \cdot x_0, \frac{1}{N}\right) + o_{\e_j}(1) \\
&=& \sum_{k=0}^{N-1}  \phi^{\e_j}\left( x_0, \left((\g^{(j,N)}_k)^{-1} \g^{(j,N)}_{k+1}\right) \cdot x_0, \frac{1}{N} \right) + o_{\e_j}(1)\\
&\geq& \sum_{k=0}^{N-1}  \widehat{\phi}^{\e_j}\left( \bp(x_0), \bp\left( \left((\g^{(j,N)}_k)^{-1} \g^{(j,N)}_{k+1}\right) \cdot x_0\right), \frac{1}{N}\right) + o_{\e_j}(1),
\end{eqnarray*}
where, as usual, $\bp: X \rightarrow \overline{ X\slash\Gamma}$ denotes the projection on the abelian cover and  $\widehat{\phi}$ the Ma\~n\'e potential associated to the projected Lagrangian.
Taking the limit as $\e_j$ goes to zero ($N$ is fixed) we obtain:
\begin{eqnarray*}
 \bb^-(\bx,1) &=& \lim_{\e_j\rightarrow 0} \phi^{\e_j}(x_0, h_{\e_j}(\bx)\cdot x_0, 1) 
\geq \sum_{k=0}^{N-1}  \frac{1}{N} \beta( \opi \log_\infty(\delta_N( (\bar{x}^{N}_k)^{-1} \bar{x}^{N}_{k+1} ) )\\
 &=& \int_0^1 \beta(\opi(\dot{\xi}_N(s)))ds,
 \end{eqnarray*}
where  $\xi_N:[0,1] \longrightarrow G_\infty$ denotes the horizontal path whose derivative is piecewise constant and equal to 
$v^N_k:= \opi \log_\infty(\delta_N( (\bar{x}^{N}_k)^{-1} \bar{x}^{N}_{k+1} ) \in V_1$ on each interval $[k/N, (k+1)/N]$ for $k=0,\ldots, N-1$ (pay attention that, although $\xi_N(0)=e$, it is not necessarily true anymore that $\xi_N(1)=\bx$). \\
Observe that since the  norm of the velocities of the minimizers $\s_j$ are equibounded (let us denote by $A$ a possible upper-bound), then 
$d_{\e_j}(\g^{(j,N)}_k \cdot x_0,  \g^{(j,N)}_{k+1}\cdot x_0) \leq \frac{A}{N} + o_{\e_j}(1)$ and therefore
$\e_j d_{\Gamma}(\g^{(j,N)}_k ,  \g^{(j,N)}_{k+1}) \leq \frac{B}{N} + o_{\e_j}(1)$, for some positive constant $B$ (which may differ from $A$, if $d_\G$ is not the orbit metric). Hence,
 $d_{\infty}(e,  \delta_N((\bar{x}^{N}_k)^{-1} \bar{x}^{N}_{k+1})) \leq B$ for $k=0,\ldots, N$. \\
 We can then apply 
\cite[Lemma 1.5 and Lemma 3.2]{BL} and deduce that $d_{\infty}(\xi_N(1), \bx) = O(1/N)$.
Hence:
\begin{eqnarray*}
\bb^-(\bx,1)  \geq \int_0^1 \beta(\opi(\dot{\xi}_N(s)))ds \geq \inf_{\s \in \cH_{\bx}} \int_0^1 \beta(\opi(\dot{\s}(s)))\,ds + O(1/N).
\end{eqnarray*}
The claim follows by taking the limit as $N\rightarrow +\infty$. 
And this concludes  the proof.\\
 \end{Proof}

\subsection{Convergence of Rescaled Solutions and Proof of the Main Theorem}\label{secproofmain}

We want to show now that the rescaled solutions $v^\e(\cdot, \cdot)$ of (\ref{HJeps}) -- see formula (\ref{vepsilon}) -- converge uniformly on compact sets  (in the sense of Definition \ref{defconv}) to 
\begin{eqnarray*}
\bar{v}(\bx,T) &=& \inf_{\by\in G_\infty} \left\{
\bbf(\by) +  \bb\left(\by^{-1}\bx, T
\right)
\right\}  
= \inf_{\by\in G_\infty} \left\{
\bbf(\by) +  T \overline{L}\left(\
\delta_{1/T}(\by^{-1}\bx)
\right)
\right\},
\end{eqnarray*}
where the second equality follows from Proposition \ref{riscalamentobeta}.\\

\begin{Rem} \label{remBCP}
It has been proven in \cite[Theorem 3.4]{BCP}  that  such a $\bar{v}$ is the unique viscosity solution to the following problem:
$$
\left\{ \begin{array}{ll}
\partial_t \bv (\bx,t) + \overline{H}(\nabla_\cH \bv(\bx,t)) = 0 & (\bx,t) \in G_\infty \times (0,T)\\
\bv(\bx,0) = \bbf(\bx) & \bx\in G_\infty,
\end{array}
\right.
$$
where  $\overline{H}: \frac{\cg_{\infty}}{[\cg_{\infty},\cg_{\infty}]} \longrightarrow  \R$ is the convex conjugate of $\b$ restricted to the subspace $\opi(\cg_{\infty}) \simeq H_1( X\slash\Gamma; \R)$ (observe that without assuming $H_1(X;\R)=0$, $\opi(\cg_{\infty})$ might be a proper subspace of $H_1( X\slash\Gamma; \R)$).  In particular, modulo identifying $H^1( X\slash\Gamma; \R)$ with  $\frac{\cG_\infty}{[\cG_\infty,\cG_\infty]}$, we can conclude that $\overline{H}$ coincides with Mather' $\alpha$-function, analogously to what happens  in \cite{LPV, CIS}.\\
\end{Rem}

%\vspace{20pt}

Let us start by proving the following Lemma (see also \cite[Lemma 3.3]{CIS}).

\begin{Lem}\label{lemmaCIS}
Given a compact set $K\subset G_\infty$ and a compact interval $J\subset (0,+\infty)$, there exists a positive constant $C>0$ such that for sufficiently small $\e>0$:
$$
d_\e(h_\e(\bx)\cdot x_0, y_\e) \leq C
$$
for any $\bx\in K$, $T\in J$ and $y_\e\in X$ realizing the minimum in the definition of $v^\e(h_\e(\bx)\cdot x_0, T)$.\\
\end{Lem}

\begin{Proof}
Since $\bbf: G_\infty \rightarrow \R$ is {assumed to have at most linear growth}, \ie
$$
\sup_{\bx \in G_\infty} \frac{|\bbf(\bx)|}{1+d_\infty(e,\bx)} <\infty,
$$
then there exist $A, B >0$ such that:
$$
f(\bx) \geq -A d_{\infty}(e, \bx) - B.
$$
Since $f_\e$ are assumed to converge locally uniformly to $\bbf$ (see Definition \ref{defconv}), if  $\bx$ is in a compact set $K\subset G_\infty$, then for $\e$ sufficiently small we have that
$f_\e(h_\e(\bx)\cdot x_0) \leq Q$, for some $Q>0$, and
\begin{eqnarray}\label{superlinfeps}
f_\e(h_\e(\bx)\cdot x_0) &\geq& \bbf(\bx) - 1
\geq -A \,d_{\infty}(e, \bx) - (B+1).
\end{eqnarray}
%where in the last step we have used Remark \ref{propertyfe}.
Recall now that $L$ is superlinear, therefore for any $M>0$ there exists $N=N(M)>0$ such that, taken any pair $x, y \in X$ and a curve $\g$ linking them in time $T$, then:
\begin{equation}\label{superlinLag}
\int_0^T  L(\g, \e\dot{\g}) dt \geq M\, d_\e(x,y) - NT.
\end{equation}
Moreover, from the definition of $v^\e$:
$$
v^\e(h_\e(\bx)\cdot x_0,T):= \inf  \left\{
f_\e(\g(0)) + \int_0^T L_\e (\g(t),\dot{\g}(t))\,dt\; \big| \; \g\in C^1([0,T], X),\; \g(T)=h_\e(\bx)\cdot x_0 \right\}.
$$
Then, if follows easily that for $\e$ small, if $\bx\in K$ and $T\in J$ then (let us denote $M_J:= \max_{T\in J}T$):
\begin{eqnarray}\label{estimateone}
v^\e(h_\e(\bx)\cdot x_0, T) \leq Q + M_J\, \max_{y\in X} |L(y, 0)| \leq Q_1.
\end{eqnarray}
Let $y_\e$ be a point achieving the infimum in $v^\e(h_\e(\bx)\cdot x_0, T)$ for some $\bx\in K$ and $T\in J$. Then, using  
that $f_\e$ {are equiLipschitz}  (with constant $\Lambda$),  the estimate in (\ref{superlinLag}) with $M$ sufficiently large (it suffices to take $M>\Lambda$) and (\ref{superlinfeps}), we obtain:
\begin{eqnarray*}
v^\e(h_\e(\bx)\cdot x_0, T) &\geq& f_\e(y_\e) + M\, d_\e(h_\e(\bx)\cdot x_0,y_\e) - NT \\
&\geq&    f_\e(h_\e(\bx)\cdot x_0)  +(M - \Lambda) d_\e(h_\e(\bx)\cdot x_0,y_\e) - NT \\
&\geq&   - A d_\infty (e, \bx) - (B+1) + (M - \Lambda) d_\e(h_\e(\bx)\cdot x_0,y_\e) - N M_J \\
&\geq&   M_1 d_\e(h_\e(\bx)\cdot x_0,y_\e) - N_1,
\end{eqnarray*}
where $M_1:= M-\Lambda$ and $N_1:= A \,{\rm diam_{\infty}}(K) + B +1 + NM_J$.
Using this estimate and (\ref{estimateone}), we can conclude that:
\begin{eqnarray*}
d_\e(h_\e(\bx)\cdot x_0,y_\e) &\leq& \frac{1}{M_1} \left( A_1+ Q_1 \right) =: C=C(J, K).
\end{eqnarray*}
\end{Proof}

%\vspace{20 pt}

\subsubsection{Proof of the Main Theorem}

We can now prove the Main Theorem (see statement in subsection \ref{secmainres}). Observe that once we have proved the convergence of $v_\e$ to $\bar{v}$ as $\e$ goes to zero (part i)), then the proof will be essentially complete: part ii), in fact, follows from Proposition \ref{whatislimit}, while part iii) follows from Remark \ref{remBCP}.\\

\begin{Proof}[Main Theorem]
Let us consider a point $(\bx, T)$ in $G_\infty\times (0,+\infty)$;  the sequence $\tilde{\delta}_\e(h_\e(\bx))$ converges to $\bx$ with respect to the metric $Q$ (see Remark \ref{propertyfe} {ii})). \\

\noindent $\bullet$ Let us first prove the pointwise convergence, {\it i.e.},
$$
\lim_{\e\rightarrow 0^+} v^\e(h_\e(\bx)\cdot x_0, T) = \bv(\bx, T).
$$
 Assume that $v^{\e_k}(h_{\e_k}(\bx)\cdot x_0, T)$ is a subsequence converging to the liminf  and let us denote by $y_{\e_k}$ the corresponding points at which the infimum (minimum) in the definition of $v^{\e_k}(h_{\e_k}(\bx)\cdot x_0, T)$ is achieved. Let us denote 
by $\g_{\e_k}$ a sequence of elements of $\Gamma$ such that $\e_k d(\g_{\e_k} \cdot x_0, y_{\e_k}) \rightarrow 0$ as $\e_k$ goes to zero (this is possible, since the set $\G\cdot x_0$ becomes denser and denser, with respect to the distance $d_{\e_k}$).
According to Lemma \ref{lemmaCIS}, $d_{\e_k}(h_{\e_k}(\bx)\cdot x_0, y_{\e_k})<C$ for $\e_k$ small, so this means that (at least for $\e_k$ small) $\tilde{\delta}_{\e_k}(\gamma_{\e_k})$ lie in a compact set (see Remark \ref{propertyfe}  ii)) and therefore -- up to possibly extracting a subsequence -- we can assume that $\tilde{\delta}_{\e_k}\g_{\e_k}$ converges to some $\by \in G_\infty$ (in the sense of Remark \ref{propertyfe} ii)).

In particular, $\tilde{\delta}_{\e_k}\big(\g_{\e_k}^{-1}h_{\e_k}(\bx)\big)$ converges to $\by^{-1}\bx$ and
\begin{eqnarray*}
&& \lim_{\e_k\rightarrow 0^+}  \phi^{\e_k}(\g_{\e_k}\cdot x_0, h_{\e_k}(\bx) \cdot x_0, T ) \;=\;
\lim_{\e_k\rightarrow 0^+}  \phi^{\e_k}(x_0, \left(\g_{\e_k}^{-1} h_{\e_k}(\bx)\right) \cdot x_0, T )\\
&& \qquad =\; \lim_{\e_k\rightarrow 0^+}  \phi^{\e_k}(x_0, h_{\e_k}(\by^{-1}\bx) \cdot x_0, T )
\;=\;   \bb(\by^{-1}\bx, T).\\
\end{eqnarray*}

\noindent Moreover,
\begin{eqnarray*}
%\lim_{\e_k\rightarrow 0} 
\left| f_{\e_k}(y_{\e_k}) - f_{\e_k}(\g_{\e_k}\cdot x_0)  \right |  & \leq &  \e_k K_1 d( y_{\e_k}, \g_{\e_k}\cdot x_0 ) \longrightarrow 0 
\end{eqnarray*}
and
\begin{eqnarray*}
%\lim_{\e_k\rightarrow 0} 
\left| \phi^{\e_k}(y_{\e_k}, h_{\e_k}(\bx)\cdot x_0, T ) - 
\phi^{\e_k}(\g_{\e_k}\cdot x_0, h_{\e_k}(\bx)\cdot x_0, T )
  \right |  & \leq &  \e_k K_2 d( y_{\e_k}, \g_{\e_k}\cdot x_0 ) \longrightarrow 0. \\
\end{eqnarray*}

\noindent Therefore,
\begin{eqnarray*}
&& \lim_{\e_k\rightarrow 0} v^{\e_k}(h_{\e_k}(\bx)\cdot x_0) \;=\;
\lim_{\e_k\rightarrow 0}
\left(
f_{\e_k}(y_{\e_k}) +  \phi^{\e_k}(y_{\e_k}, h_{\e_k}(\bx)\cdot x_0, T )\right) \\
&& =\;
\lim_{\e_k\rightarrow 0}
\left(
f_{\e_k}(\g_{\e_k}\cdot x_0) +  \phi^{\e_k}(\g_{\e_k}\cdot x_0, h_{\e_k}(\bx)\cdot x_0, T )\right) \;
=\;
\bbf(\by) + \lim_{\e_k\rightarrow 0} \phi^{\e_k}(\g_{\e_k}\cdot x_0, h_{\e_k}(\bx)\cdot x_0, T ) \\
&&=\; 
\bbf(\by) + \bb(\by^{-1}\bx, T) \geq \bv(\bx, T).
\end{eqnarray*}
Therefore, 
\begin{equation}\label{convliminf}
\liminf_{\e\rightarrow 0^+} v^\e(h_\e(\bx)\cdot x_0) \geq  \bv(\bx, T).\\
\end{equation}

\vspace{10 pt}

Let now $\by$ be the optimal element in the definition of $\bv(\cdot, T)$  (it does exist since $\bL$ is convex and superlinear) and  consider the sequence
$h_\e(\by)$.%, which is such that  $\delta_\e(h_\e(\by))$ converges to $\by$. 
Then:
\begin{eqnarray*}
v^\e(h_\e(\bx)\cdot x_0, T) &\leq& f_\e(h_\e(\by)\cdot x_0) + \phi^\e(h_\e(\by)\cdot x_0,h_\e(\bx)\cdot x_0, T). 
\end{eqnarray*}
Taking the limit as $\e$ goes to zero we obtain:
\begin{eqnarray*}
&& \limsup_{\e \rightarrow 0^+} v^\e(h_\e(\bx)\cdot x_0, T) \;\leq\; \limsup_{\e\rightarrow 0^+} \left(f_\e(h_\e(\by)\cdot x_0) + \phi^\e(h_\e(\by)\cdot x_0,h_\e(\bx)\cdot x_0, T) \right)\\
&& \quad \leq\; \bar{f}(\by) + 
\lim_{\e\rightarrow 0^+} \phi^\e(x_0, (h_\e(\by)^{-1} h_\e(\bx))\cdot x_0, T)
\;=\; \bar{f}(\by) + \bb(\by^{-1}\bx, T) = \bv(\bx, T).
\end{eqnarray*}

This and (\ref{convliminf}) allow us to conclude that for each $T>0$ and $\bx\in G_\infty$ we have
$$
\lim_{\e \rightarrow 0^+} v^\e(h_\e(\bx)\cdot x_0, T) = \bv (\bx, T).
$$

%\vspace{20 pt}

\noindent $\bullet$ Let now  $\bx_\e \rightarrow \bx$ in $(G_\infty, d_\infty)$ and $T_\e\rightarrow T$. Then, using the equiLipschitzianity of the $v^\e$ (see Proposition \ref{equilip}, observing that for $T$ in a compact set of $(0,+\infty)$ the Lipschitz constant can be chosen uniformly, at least for small $\e$):
\begin{eqnarray*}
&& |v^\e(h_\e(\bx_\e)\cdot x_0, T_\e) -\bv(\bx, T)| \;\leq\;
|v^\e(h_\e(\bx_\e)\cdot x_0, T_\e) - v^\e(h_\e(\bx)\cdot x_0, T_\e)| \\
&& \qquad  \;+\; |v^\e(h_\e(\bx)\cdot x_0, T_\e) - v^\e(h_\e(\bx)\cdot x_0, T)|  \;
+\;  |v^\e(h_\e(\bx)\cdot x_0, T) -\bv(\bx, T)| \\
&& \qquad \leq\;  K_1 d_\e ( h_\e(\bx_\e)\cdot x_0, h_\e(\bx)\cdot x_0) + K_2 |T_\e - T| 
\;+\; |v^\e(h_\e(\bx)\cdot x_0, T) -\bv(\bx, T)| \\
&& \qquad \leq\;  K_3 d_\infty (\bx_\e, \bx) + K_2 |T_\e - T| + |v^\e(h_\e(\bx)\cdot x_0, T) -\bv(\bx, T)|,
\end{eqnarray*}
which goes to zero as $\e$ goes to zero. Hence:
\begin{equation}\label{limite}
\lim_{\e \rightarrow 0^+} v^\e(h_\e(\bx_\e)\cdot x_0, T_\e) = \bv (\bx, T).\\
\end{equation}

%\vspace{20 pt}

\noindent $\bullet$ Let us now prove that this convergence is uniform for $\bx$ in compact sets of $G_\infty$ and $T$ in a compact interval $J\subset (0,+\infty)$, \ie
$$
\lim_{\e\rightarrow 0^+} \sup_{B_\infty(R)\times J} |v^\e(h_\e(\bx)\cdot x_0, T)) - \bv (\bx, T)| =0.
$$

If by contraddiction this were not true, then we could find a sequence $\e_k\rightarrow 0$,  elements $(\bx_k, T_k ) \in B_\infty(R)\times J$ and $\delta>0$, such that
\begin{equation}\label{assurdo}
|v^{\e_k}(h_{\e_k}(\bx_k)\cdot x_0, T_k) - \bv (\bx_k, T_k)| >\delta \qquad \forall \;k.\\
\end{equation}

Up to extracting a subsequence, we could assume that $x_k \rightarrow \bx \in B_\infty(R)$ and $T_k\rightarrow T>0$, and it would follow from (\ref{limite}) that
$$
\lim_{\e_k \rightarrow 0} v^{\e_k}(h_{\e_k}(\bx_k), T_k) = \bv (\bx, T),
$$
which clearly contradicts (\ref{assurdo}). Therefore, the convergence must be uniform.\\
\end{Proof}

\appendix

\section{The Effective Hamiltonian and Mather--Ma\~n\'e Theory} \label{appMather}
It turns out that the effective Hamiltonian -- that we have defined in subsection \ref{ss103} in terms of  the solutions to the cell problem -- is also extremely significant from a dynamical systems point of view, particularly in the study of the associated Hamiltonian dynamics by means of variational methods: what is nowadays known as Mather and Ma\~n\'e theory. We refer interested readers to \cite{Mather91, Mane, Gonzalobook, SorrentinoLectureNotes} for more detailed presentations of these topics.

Roughly speaking, this theory is based on the study of  particular orbits and invariant measures of the flow that are obtained as minimizing solutions to variational problems related to the so-called {\it Principle of least (Lagrangian) action}. As a result of this, these objects present a much richer structure and rigidity than one might generally expect, and the corresponding invariant sets -- the so-called Mather, Aubry and Ma\~n\'e sets --  play an important  role in determining both the local and the global dynamics of the system.

In this setting, the value of the effective Hamiltonian appears in many noteworthy forms and has consequently been named in different ways by the various communities: {\it minimal average action, Mather's $\alpha$-function, Ma\~n\'e critical values}, etc...\\

Let us briefly recall some of these definitions.
In what follows, let  $L:TM\longrightarrow \R$ denote the Lagrangian associated to $H$ (given by Legendre-Fenchel duality)  and  for any cohomology class $c\in H^1(M;\R)$ let us consider
a closed $1$-form $\eta_c$ representing it (it is easy to check that the definitions below do  not depend on the chosen representative, but only on their cohomology class). Then:\\

\begin{itemize}
\item[{\small 1)}] If $\calM_{L}$ denotes the set of invariant probability measures for the Euler-Lagrange flow associated to  $L$, then:
\beqano
\bH(c) = - \min_{\m \in \calM_L} \int_{\rT M} (L(x,v) -\eta_c(x)\cdot v)\, d\m.
\eeqano
The value of the right-hand side is usually denoted by $\alpha(c)$ and the collection of these values $\a: H^1(M;\R) \longrightarrow  \R$
is what is known as {\it Mather's $\alpha$-function} or {\it Mather's minimal average action}.\\

\item[{\small 2)}] For any absolutely continuous curve $\g:[a,b]\longrightarrow M$, we define its {\it $L$-action} as
$$A_{L,\eta_c}(\g):= \int_a^b \Big( L(\g(t),\dot{\g}(t)) - \eta_c(\gamma(t))\cdot \dot{\gamma}(t) \Big)\,dt.$$ 
Then:
\beqano
\bH(c) &=&\inf  \{k\in\R:\; A_{L+k,{\eta_c}}(\g)\geq 0, \quad \forall\;\mbox{abs. cont. loop}\; \g\} \nonumber\\
&=&\sup\{k\in\R:\;  A_{L+k,{\eta_c}}(\g)< 0 \;  \mbox{for some abs. cont. loop}\; \g \}\,. 
\eeqano
The values on the right-hand sides are often called {\it Ma\~n\'e critical values}.\\

\item[{\small 3)}] In \cite{Carneiro} Dias Carneiro proved that $\bH(c)$ represents the energy (\ie the value of the Hamiltonian) of  action-minimizing measures or action-minimizing orbits of cohomology class $c$.\\

\item[{\small 4)}] It was proved in \cite{PPS} that  $\bH(c)$ represents the infimum of the energy values $k$'s such that the energy sublevel $\{H(x,p)\leq k\}$ contains in its interior a smooth Lagrangian graph of cohomology class $c$. In particular, it corresponds to the smallest energy sublevel containing Lipschitz Lagrangian graphs of cohomology class $c$.\\

\item[{\small 5)}] $\a: H^1(M;\R) \longrightarrow  \R$ is a convex function, so one can consider its {\it Fenchel-Legendre conjugate} defined on the dual space $(H^1(M;\R))^* \simeq H_1(M;\R)$, namely the {\it first homology group of the manifold}:
\begin{eqnarray*}
\beta: H_1(M;\R) &\longrightarrow& \R \\
h &\longmapsto & \beta(h)= \max_{c\in H^1(M;\R)} \left( \langle c,h \rangle - \alpha(c) \right),
\end{eqnarray*}
where $\langle \cdot, \cdot \rangle$ denotes the pairing between $H^1(M;\R)$ and $H_1(M;\R)$.
It turns out that also this function has a dynamical meaning; namely for every $h\in H_1(M;\R)$ it represents the minimal $L$-action of all invariant probability measures with {\it rotation vector} (or {\it Schwartzmann asymptotic cycle}) equals to $h$. We refer the reader to 
\cite{Mather91, SorrentinoLectureNotes} for more details. \\

\end{itemize}

\begin{Rem}
Finally, it is also interesting to observe that the homogenized Hamiltonian $\bH$ coincides with the {\it symplectic homogenized Hamiltonian} defined by Viterbo in \cite{Viterbo} for Hamiltonians on $T^*\T^n$. This definition was later extended to Hamiltonian on general compact manifolds in \cite{MVZ}.\\
\end{Rem}

%%%%%%%%%%%%%%%%%%%%%%%%%%%%%%%%

\vspace{1.truecm}

\end{document}